\documentclass[12pt]{article}
\usepackage{latexsym,amsmath}
\usepackage[authoryear]{natbib}
       
\def\R{{\rm I \mkern-2.5mu \nonscript\mkern-.5mu R}}
\def\N{{\rm I \mkern-2.5mu \nonscript\mkern-.5mu N}}

\newcommand{\be}{\begin{equation}}
\newcommand{\ee}{\end{equation}}
\newcommand{\bea}{\begin{eqnarray}}
\newcommand{\eea}{\end{eqnarray}}
\newcommand{\bean}{\begin{eqnarray*}}
\newcommand{\eean}{\end{eqnarray*}}

\newcommand{\qed}{$\hfill\Box$}

\newcommand{\sect}[1]{\section{#1} \setcounter{equation}{0}
                      \setcounter{table}{0} \setcounter{figure}{0}}

\newtheorem{defi}{Definition}[section]
\newtheorem{ex}[defi]{Example}
\newtheorem{con}[defi]{Condition}
\newtheorem{lem}[defi]{Lemma}
\newtheorem{theo}[defi]{Theorem}

\title{Efficient estimation for ergodic diffusion processes sampled 
at high frequency}

\author{Michael S\o rensen \\ 
{\small Department of Mathematical Sciences} \\
{\small University of Copenhagen} \\
{\small Universitetsparken 5,
DK-2100 Copenhagen {\O}, 
Denmark}}

\begin{document}

\maketitle

\begin{abstract}
A general theory of efficient estimation for ergodic diffusion processes sampled 
at high frequency with an infinite time horizon is presented. High
frequency sampling is common in  many applications, with finance
as a prominent example. The theory is formulated  
in term of approximate martingale estimating functions and
covers a large class of estimators including most of the previously 
proposed estimators for diffusion processes. Easily checked conditions
ensuring that an estimating function is an approximate martingale are
derived, and general conditions ensuring consistency and asymptotic
normality of estimators are given. Most importantly, simple conditions
are given that ensure rate optimality and efficiency. Rate optimal
estimators of parameters in the diffusion coefficient converge faster
than  estimators of drift coefficient parameters because they take
advantage of the information in the quadratic variation. The conditions 
facilitate the choice among the multitude of estimators that have been 
proposed for diffusion models. Optimal martingale estimating functions in the
sense of Godambe and Heyde and their high frequency approximations
are, under weak conditions, shown to satisfy the conditions for rate
optimality and efficiency. This provides a natural feasible method of
constructing explicit rate optimal and efficient estimating functions
by solving a linear equation.     

\end{abstract}

\newpage

\sect{Introduction}

Dynamic phenomena affected by random noise are often modelled
in continuous time by stochastic differential equations. Among the 
advantages of this approach are interpretable model parameters and
easy communication with other scientists by using a  common modelling
tool, viz.\ differential equations. A few examples are applications in
the areas of animal movement (\cite{animal}), 
climate research (\cite{2ditlevsen}), finance (\cite{ckls}, \cite{crypto}), protein
structure evolution (\cite{protein}), neuro science (\cite{bibbona2010},
\cite{FitzHugh}, \cite{biomath}), transmission of infectious diseases
(\cite{GuyLaredoVergu}, \cite{covid}) and physiology
(\cite{picchini2008}). While the dynamics is formulated in continuous
time, observations are made at discrete points in time. This
complicates statistical inference for these models, which is an intensive 
area of research, where a profusion of estimators have been proposed for
parametric diffusion models, see e.g.\ \cite{MS12} and \cite{yuima}. 
Many simulation studies have been performed to compare the relative
merits of estimators, but have not provided a clear general
picture. The simple and easily checked criteria for efficiency and
rate optimality of estimators obtained in this paper are useful for
identifying the best estimators and explain findings in simulation studies.

We consider a scalar diffusion given by the stochastic differential equation
\be
\label{sde}
dX_t=b(X_t;\alpha)dt+\sigma(X_t;\beta)dW_t,
\ee
where $(\alpha, \beta) = \theta \in \Theta \subseteq \R^2$ are
parameters to be estimated. The restrictions to a scalar process and
to two scalar parameters are made to simplify the presentation. The
results can be generalized to multivariate diffusions
as indicated in Section \ref{efficientsection}, and to multivariate
parameters by considering estimating functions of the same
dimension as the parameter vector and replacing partial derivatives by
vectors or matrices of partial derivatives. 
The process $X$ is assumed to be observed at equidistant time points $i \Delta_n$, 
$i = 0, \ldots , n$, and we consider the high frequency/infinite time horizon 
asymptotic scenario, where
\[
n \rightarrow \infty, \hspace{6mm} \Delta_n \rightarrow 0, \hspace{6mm} 
n \Delta_n \rightarrow \infty.
\]
The length of the time interval in which observations are 
made goes to infinity, which is necessary to ensure that the drift
parameter $\alpha$ can be estimated consistently. At the same time
the sampling frequency goes to infinity. This is particularly
important for diffusion processes, because the quadratic variation of a
diffusion process contains information about the parameter $\beta$ in the
diffusion coefficient, which good estimators can use when the sampling
frequency is sufficiently high. For such estimators, $\beta$ is
estimated at a higher rate than the $1/\sqrt{n\Delta_n}$-rate, which is
optimal for the drift parameter $\alpha$; see
\cite{gobet} where it is shown that the model (\ref{sde}) is locally
asymptotically normal under the high frequency/infinite time horizon 
asymptotic scenario.  In the present paper we give easily checked
conditions for rate optimality and efficiency. If the drift
coefficient is known, consistent estimators can be found without the
infinite time horizon assumption. Results on rate optimality and efficiency
when the sampling interval is bounded are given in \cite{nmms1}. 
High frequency asymptotics is often relevant in applications, because
the sampling frequency needs not be particularly high for the
asymptotic optimality results to work. It only needs to be high
relative to the characteristic time scale of the diffusion process.
For some types of economic  data, even weekly observations
can in be considered a high sampling frequency, see e.g.\ \cite{larsen}.

Our theory is phrased in terms of estimating functions of the general form
\be
\label{estfct}
G_n(\theta) = \sum_{i=1}^n g(\Delta_n, X_{i \Delta_n},
X_{(i-1)\Delta_n }; \theta),
\ee
where the function $g(\Delta,y,x; \theta)$, with values in $\R^2$, is
such that $G_n$ is, exactly or approximately, a martingale estimating
function. Specifically, $E_\theta(g(\Delta, X_{i \Delta_n},
X_{(i-1)\Delta_n }; \theta) \, | \, X_{(i-1)\Delta_n })$ is 
equal to zero or of order $\Delta^\kappa$ for some $\kappa \geq
2$. Estimators are obtained as solutions to the estimating equation
$G_n(\theta)=0$. We call such an estimator a $G_n$-{\it estimator}. For estimating
functions that are not exact martingales, the extra condition that $n
\Delta^{2(\kappa-1)} \rightarrow 0$ is needed to ensure the asymptotic results.

The theory developed here covers a large class of estimators for diffusion 
processes including most of the previously proposed estimators. The 
few that are not covered are likely to be less efficient, because non-martingale
estimating functions, in general, do not approximate the score
function as well as martingales. In particular, the theory 
covers the martingale estimating functions proposed 
by \cite{bmb&ms1} and \cite{kessler&ms}, GMM-estimators based on 
conditional moments, \cite{hansen82, hansen2, hansen93}, and the
maximum likelihood estimator and Bayesian estimators; for numerical
methods to calculate these estimators, see \cite{pedersen},
\cite{roberts}, \cite{ait-sahalia}, \cite{durhamgallant}, \cite{ait-sahalia-mykland0},
\cite{beskos9}, \cite{golightly&wilkinson2}, \cite{bladtsorensen},
\cite{MeulenSchauer} and \cite{bladtmidersorensen}. The  
pseudo-likelihood function obtained from the Gaussian Euler approximation 
to the transition density is covered too. Estimators closely related to 
the Euler pseudo-likelihood were considered by \cite{zmirou},
\cite{nakahiro3} and \cite{uchida10}. These and pseudo-likelihood functions
based on more accurate Gaussian approximations to the likelihood
function, such as those considered by \cite{kessler97},
\cite{uchidayoshida12} and \cite{uchidakitagawa14}, are also covered. \cite{uchida} and
\cite{agms} considered the Euler pseudo likelihood under a combination of high 
frequency and small diffusion asymptotics, where the diffusion 
coefficient goes to zero as $n \rightarrow \infty$. The latter
condition replaces the infinite time horizon condition.

The following condition on the function $g(\Delta,y,x; \theta)$
ensures {\it rate optimality} of estimators.
\begin{con}
\label{con3}
$ \ $
\vspace{0.5mm}
\be
\label{jacobsen}
\partial_y g_2(0,x,x;\theta) = 0
\ee
for all $x$ in the state-space of the diffusion process and all
$\theta \in \Theta$. 
\end{con}
By $\partial_y g_2(0,x,x;\theta)$ we mean $\partial_y
g_2(0,y,x;\theta)$ evaluated at $y=x$. Here $\partial_y$ denotes the
partial derivative w.r.t.\ $y$, and $g_i$ is the $i$th
coordinate of $g$.  More precisely, it is sufficient that a linear combination of the two
coordinates of $g$ satisfies (\ref{jacobsen}), but without loss of
generality it can be assumed to be $g_2$. This will be explained in
Section \ref{model}.  We will refer to (\ref{jacobsen})
as Jacobsen's condition because it equals one of the 
conditions for small $\Delta$-optimality in the sense of \cite{mj} of
martingale estimating functions; see \cite{mj1}. Jacobsen considered an
asymptotic scenario where the time between observations $\Delta$ does
not depend on $n$. In his approach the
Condition \ref{con3} was introduced to avoid a singularity in the
asymptotic variance of the estimators at $\Delta = 0$.
In our high frequency approach, the condition implies rate optimality
for estimation of the diffusion coefficient parameter.

\newpage 

Our condition for {\it efficiency} is
\begin{con}
\label{conefficiency}
$ \ $
\vspace{0.5mm}
\be
\label{drifteff}
\partial_y g_1 (0,x,x;\theta) = \partial_\alpha b(x;\alpha)/\sigma^2(x;\beta)  
\ee
and
\be
\label{diffeff}
\partial^2_y g_2 (0,x,x;\theta) = \partial_\beta
\sigma^2(x;\beta)/\sigma^4(x;\beta),
\ee
for all $x$ in the state space of the diffusion process and all
$\theta \in \Theta$.
\end{con}

The Conditions \ref{con3} and \ref{conefficiency} are, under weak
regularity conditions, shown to be satisfied by 
martingale estimating functions that are optimal in the sense of
\cite{godambeheyde}, which is both very useful and quite 
surprising. Useful because it provides an easy method of constructing
rate optimal and efficient estimators, and surprising because Godambe-Heyde
optimality is a local property in the sense that it is a property of a
particular class of estimating functions. Therefore, there is no
a priori reason to except this property to imply global properties
like rate optimality and efficiency.  
Martingale estimating functions give consistent estimators at all 
sampling frequencies, see \cite{bibbyjacobsenms}, 
and Godambe-Heyde optimal martingale estimating functions 
are known to often provide estimators with a high efficiency, see e.g.\
the simulation studies in \cite{overbeck} and \cite{larsen}.
The results in this paper explain why this is the case. 

The paper is organized as follows. Section 2 sets up the model, the
class of approximate martingale estimating functions, and the
assumptions and the notation used throughout the paper. A number of
well-known estimators are shown to be covered by the theory, and a
lemma of independent interest gives fundamental identities and
characterizes approximate martingale estimating functions of order
$\kappa$. Section 3 develops the high frequency asymptotic 
theory for general estimating functions as well as for estimating
functions satisfying Condition \ref{con3}. The conditions for
efficiency are derived in Section 4, and it is proved
that Godambe-Heyde optimal martingale estimating functions and their
high frequency approximations are rate optimal and efficient, which
provides a feasible method of constructing explicit rate optimal and
efficient estimating functions. Examples are considered, including 
the Euler pseudo-likelihood and maximum likelihood estimation. Proofs
and some lemmas are given in Section 5, where also tools for studying
high frequency asymptotic properties of estimators are
provided. Section 6 concludes.

\sect{Model, conditions and notation}
\label{model}

We consider observations $X_{t_0^n}, \ldots, X_{t_n^n}$ of the process 
given by (\ref{sde}) at the time points $t_i^n =i \Delta_n$, $i = 0, 
\ldots , n$. We suppose that a solution of the stochastic
differential equation (\ref{sde}) exists, is unique in law, and
is adapted to the filtration generated by the Wiener process 
$W$ and the initial value $X_0$. The state-space
of $X$ is denoted by $(\ell, r)$, where $-\infty \leq \ell < r \leq
\infty$, and we assume that $v(x;\beta) = \sigma^2(x;\beta)>0$ for all
$x \in (\ell, r)$. Furthermore, we assume that $\theta =
(\alpha,\beta) \in \Theta$, where $\Theta$ is a subset of $\R^2$, and
that the true parameter value $\theta_0 = (\alpha_0,\beta_0) \in
\mbox{int} \, \Theta$, the interior of $\Theta$. It is no
serious restriction to assume that $\Theta$ is convex.

A function $f(y,x;\theta)$ is said to be of {\it polynomial growth} 
in $y$ and $x$ uniformly for $\theta$ in a compact set if, for any
compact subset $K \subseteq \Theta$, there exists a constant $C > 0$ such 
that $\sup_{\theta \in K}|f(y,x;\theta)| \leq C (1+|x|^C+|y|^C)$ for
all $x, y \in (\ell, r)$. The assumptions of
polynomial growth in this paper are made to simplify the theory. These
assumptions are satisfied for most models used in practice, but could
be relaxed. 

Here and in the rest of the paper,
$R(\Delta,y,x; \theta)$ denotes a (generic) function such that
$|R(\Delta,y,x; \theta)|$ $\leq F(y,x;  \theta)$, for all $\Delta$,
where $F$ is some function of polynomial growth in $y$ and $x$ uniformly for
$\theta$ in a compact set. Similarly for $R(\Delta, x; \theta)$.

\vspace{2mm}

\begin{defi}
We let $C_{p,k_1,k_2,k_3}(\R_+ \times (\ell, r)^2 \times \Theta)$ 
denote the class of real functions $f(t,y,x;\theta)$ satisfying that
\begin{itemize}

\item[(i)] $f(t,y,x;\theta)$ is $k_1$ times continuously 
differentiable with respect to $t$, $k_2$ times continuously 
differentiable with respect to $y$, and $k_3$ times continuously differentiable 
with respect to $\theta$

\item[(ii)] $f$ and all partial derivatives $\partial^{i_1}_t \, 
\partial^{i_2}_y \, \partial^{i_3}_\alpha \, \partial^{i_4}_\beta f$,
$i_j =0, \ldots k_j$, $j=1,2$, $i_3 + i_4 \leq k_3$ are continuously differentiable 
with respect to $x$ and are of 
polynomial growth in $x$ and $y$ uniformly for $\theta$ in 
compact sets (for fixed $t \leq 1$) 

\item[(iii)] $f$ has an expansion
\be
\label{gdeltaexp}
f(\Delta,y,x; \theta) = \sum_{i=0}^{k_1} \frac{\Delta^i}{i!} f^{(i)}(y,x;\theta)
+ \Delta^{k_1+1} R(\Delta,y,x;\theta).
\ee

\end{itemize}

\vspace{1mm}

The classes $C_{p,k_2,k_3}((\ell, r) \times \Theta)$ and $C_{p,k_2,k_3}
((\ell, r)^2 \times \Theta)$ are defined similarly  (with property
(iii) omitted) for functions of the form $f(y;\theta)$ and
$f(y,x;\theta)$, respectively.   
\end{defi}

%\vspace{3mm}

We assume that the stochastic
differential equation (\ref{sde}) satisfies the following condition.
\begin{con}
\label{con1}
The following holds for all $\theta \in \Theta$:

\begin{itemize}

\item[(1)]
\be
\label{infinitescale}
\int_{x^\#}^r s(x;\theta) dx = \int^{x^\#}_{\ell} s(x;\theta) dx =
\infty
\ee
and
\be
\label{moments1}
\int_{\ell}^r x^k \tilde{\mu}_\theta (x) dx < \infty
\ee
for all $k \in \N$, where $x^\#$ is an arbitrary point in $(\ell,r)$,
\be
\label{scale}
s(x;\theta)=\exp \left( -2 \int_{x^\#}^x
  \frac{b(y;\alpha)}{v(y;\beta)} dy \right)
\ee
and
\be
\tilde{\mu}_\theta (x) = [s(x;\theta)v(x;\beta)]^{-1}
\ee

\item[(2)] $\sup_t E_\theta(|X_t|^k ) < \infty$ for all $k \in \N$

\item[(3)] $b, \sigma \in C_{p,4,1}((\ell, r) \times \Theta)$  

\item[(4)] There exists a constant $C_\theta$ such that for all $x,y
  \in (\ell,r)$
\[
| b(x;\alpha) - b(y;\alpha) | + | \sigma(x;\beta) - \sigma(y;\beta)|
\leq C_\theta |x -y|.
\]

\end{itemize}

\end{con}

\vspace{1mm}

The conditions (\ref{infinitescale}) and (\ref{moments1}) with $k=1$ 
ensure that the process $X$ is ergodic with an invariant probability
measure with Lebesgue density
\be
\label{invmeasure}
\mu_\theta(x) = \tilde{\mu}_\theta (x)/\int_\ell^r 
\tilde{\mu}_\theta (y)dy.
\ee
If $X$ is stationary, Condition \ref{con1} (2) is obviously
satisfied under (\ref{moments1}). Similarly, if $X$ is sufficiently
mixing that all moments converge as $t \rightarrow \infty$.

We consider estimating functions of the general form (\ref{estfct})
where the function $g(\Delta,y,x; \theta)$ has values in $\R^2$ and
satisfies the following condition.

\begin{con}
\label{con2}
$ \ $
\vspace{0.5mm}
\begin{itemize}

\item[(1)] There exists a $\kappa \geq 2$ such that
\be
\label{mart}
\hspace{-1mm} E_\theta (g_i(\Delta_n, X_{t_j^n}, X_{t_{j-1}^n}; 
\theta) \, | \, X_{t_{j-1}^n}) = \Delta_n^\kappa \, R(\Delta_n,X_{t_{j
   -1}^n};\theta) \ \mbox{ for } i=1,2 \mbox{ and all } \theta \in \Theta
\ee

\item[(2)] $ \ $
\vspace{-7mm}  
\bean
g_i(\Delta,y,x; \theta) &\in& C_{p,2,6,2}(\R_+ \times (\ell, r)^2 \times
\Theta), \ i=1,2 \\
g_i^{(j)}(y,x; \theta) &\in& C_{p,2(3-i),2}((\ell, r)^2 \times 
\Theta), \ i=1,2, \ j=0,1,2, 
\eean
where the $g^{(j)}$s are the functions appearing in
the expansion (\ref{gdeltaexp}).

\end{itemize}

\end{con}

We call an estimating function satisfying Condition \ref{con2} (1) an
{\it approximate martingale estimating function of order $\kappa$}.
For any non-singular $2 \times 2$ matrix, $M(\Delta, n, \theta)$, the
estimating functions $M(\Delta_n, n, \theta) G_n(\theta)$ and
$G_n(\theta)$ give identical estimators.  
We call them {\it versions} of the same estimating function. Since the matrix 
$M$ may depend on $\Delta_n$, not all versions satisfy Condition
\ref{con2} and other conditions in the paper, in particular Conditions
\ref{con3}  and \ref{conefficiency}. We say that
(\ref{estfct}) is an approximate martingale estimating function of
order $\kappa$, if there exists a version which satisfies
(\ref{mart}), and for which the limit $g(0,y,x;\theta)$ is finite
with,  for at least one value of $(x,y,\theta)$, all coordinates
different from zero. We use this version   
in the proofs. It is, typically, obtained by multiplying one or
both of the coordinates by a power of  $\Delta_n$; examples are given
in Section \ref{efficientsection}.

It could be argued, that an estimating functions that satisfies
(\ref{mart}) with $\kappa = 1$ could equally well be called an approximate
martingale. However, the asymptotic theory in this case is entirely
different from the case $\kappa \geq 2$ and requires a separate
study. Some particular examples are studied in \cite{jorgensensorensen}. 

The generator of the solution to (\ref{sde}) is the differential 
operator
\be
\label{generator}
L_\theta = b(x;\alpha) \frac{d}{dx} + \dfrac12
v(x;\beta) \frac{d^2}{dx^2}
\ee 
Here we take
the domain of $L_\theta$ to be the set of all twice continuously
differentiable functions defined on the state space.
For $f \in C_{p,2(k+1)}((\ell, r))$ and $b,\sigma \in C_{p,2k,0}((\ell, r) \times \Theta)$,
\be
\label{expansion}
\mbox{E}_\theta (f(X_{t+\Delta}) \, | \, X_t) = 
\sum_{i=0}^k \frac{\Delta^i}{i!} L_\theta^i f (X_t) + \Delta^{k+1}R(\Delta,X_t;\theta) , 
\ee
where
\[
\Delta^{k+1}R(\Delta,X_t;\theta) = \int_0^\Delta  \int_0^{u_1} \cdots 
\int_0^{u_k} E_\theta( L_\theta^{k+1} f (X_{t+u_{k+1}}) \, | \, X_t)
du_{k+1} \cdots du_1, 
\]
see e.g.\ \cite{MS12}. The properties of the remainder term follow from
Lemma \ref{lemma5} in Section \ref{proofs}. 
When we apply the generator to a function $h(y,x)$ of two variables,
we mean
\be
\label{generator2}
L_\theta (h)(y,x) = b(y;\alpha) \partial_y h(y,x)+ 
\mbox{ \small $\frac{1}{2}$} v(y;\beta) \partial^2_{y} h(y,x),
\ee
and for a function $h(\Delta,y,x;\theta)$ that depends also on 
$\Delta$ and $\theta$, we use the notation
\[
L_\theta (h(\Delta;\tilde{\theta}))(y,x) = b(y;\alpha) \partial_y h(\Delta,
y,x ;\tilde{\theta})+ \mbox{ \small $\frac{1}{2}$} v(y;\beta) 
\partial^2_{y} h(\Delta, y,x ; \tilde{\theta}).
\]

The following lemma provides identities that play an essential role in
the proofs of the asymptotic theory in the next section. Note that
$L_\theta$ is applied coordinate-wise to a vector valued function, and
that, depending on the context, $0$ can also denote a $0$-vector.

\begin{lem}
\label{lemma1}
Let $G_n$ be an estimating function of the form (\ref{estfct}), where $g_i \in
C_{p,\kappa-1,2 (\kappa-1),0}(\R \times (\ell, r)^2 \times \Theta)$, i=1,2, for a
$\kappa \geq 2$, and assume Condition \ref{con1}. 

Then $G_n$ is an approximate martingale estimating function of order
$\kappa \geq 2$ (i.e.\ it satisfies (\ref{mart})) if and only if 
\[
\sum_{i=0}^k {k \choose i}L_\theta^{k - i}(g^{(i)}(\theta))(x,x)
=0, \ \ \ k = 0, \ldots, \kappa-1,
\]
for all $x \in (\ell,r)$ and $\theta \in \Theta$ (the $g^{(i)}$s are
the functions in the expansion $(\ref{gdeltaexp})$).

In particular, if $G$ is an approximate martingale estimating
function, then
\bea
\label{lem1-1}
g^{(0)}(x,x;\theta) &=& 0 \\
\label{lem1-2}
g^{(1)}(x,x;\theta) &=& - L_\theta (g^{(0)}(\theta))(x,x)
\eea
for all $x \in (\ell,r)$ and $\theta \in \Theta$.
\end{lem}

\subsection{Examples}

The prototype of estimating functions satisfying the condition 
(\ref{mart}) are the {\it martingale estimating functions} for which
\[
E_\theta (g(\Delta_n, X_{t_i^n}, X_{t_{i-1}^n}; \theta) \, 
| \, X_{t_{i-1}^n}) = 0.
\]
They often have the form
\be
\label{specielestfct}
g(\Delta,y,x; \theta) = A(x, \Delta;\theta) \left[f(y;\theta) -
  \pi_\theta^{1,\Delta} f (x;\theta)  \right], 
\ee
where $f(y;\theta) = (f_1(y;\theta), \ldots , f_N(y;\theta))^T$, with $f_i$
real-valued, $A(x,\Delta;\theta)$
a $2 \times N$-matrix of weights, and $\pi_\theta^{1,\Delta}$ denotes the 
transition operator given by
\be
\label{transitionop}
\pi_\theta^{1,\Delta} h (x) = 
E_\theta (h(X_\Delta) \, | \, X_0 = x)
\ee
for a real-valued function $h$.
Here and later $x^T$ denotes the transpose of a vector or  matrix
$x$. The weight matrix $A$ can be chosen suitably, for instance to obtain 
rate optimality and efficiency. Examples are polynomial estimating functions, where
the real functions $f_j$ are power functions or more general polynomials.
The {\it quadratic martingale estimating function}, obtained for
$N=2$, $f_1(x)=x$ and $f_2(x)=x^2$, is as a useful simple example, see Section
\ref{efficientsection}. Polynomial estimating functions are
particularly useful for the class of Pearson diffusions, for which all (finite)
moments (conditional as well as unconditional) can be found
explicitly, see \cite{fs:08}. Other instances  
are the estimating functions based on eigenfunctions of the 
generator (\ref{generator}) proposed by \cite{kessler&ms}. 

The econometric {\it generalized method of moments} (GMM) 
based on conditional moments is covered by our
theory. This method is in practice often implemented as
follows; see \cite{clm}. The starting point
is an $N$-dimensional function $h(\Delta,y,x; \theta)$ for which each 
coordinate satisfies that $E_\theta (h_j(\Delta_n, X_{t_i^n}, 
X_{t_{i-1}^n}; \theta) \, | \, X_{t_{i-1}^n}) = 0.$ Let $A_n$ be
an $N \times N$-matrix such that $m_n(\theta) = A_n \sum_{i=1}^n
h(\Delta_n, X_{t_i^n}, X_{t_{i-1}^n}; \theta)$ converges in probability.
For the usual low frequency asymptotics, where $\Delta_n$ does not 
depend on $n$, $A_n = n^{-1} I_N$ ($I_N$ denotes the identity matrix),
but for the high frequency asymptotics considered in this paper, a
different choice of $A_n$ is usually necessary, as will 
become clear in the next section. The GMM-estimator is
obtained by minimizing $Q_n(\theta) = m_n(\theta)^T W_n \, m_n(\theta)$,
where $W_n$ is an $N \times N$-matrix such that $W_n \rightarrow W$ in
probability. It is typically the (suitably normalized) inverse
of a consistent estimator of the covariance matrix of $m_n(\theta)$.
Under weak regularity conditions, the GMM-estimator solves the estimating
equation $\partial_\theta Q_n(\theta) = \partial_\theta m_n(\theta)^T W_n \,
m_n(\theta)=0$, so if $\partial_\theta m_n(\theta) \rightarrow D(\theta)$ 
in probability (a necessary condition for asymptotic results about the
GMM-estimator), then the GMM-estimator has the same asymptotic
behavior as the estimator obtained from $D(\theta)^T W A_n
\sum_{i=1}^n h(\Delta_n, X_{t_i^n}, X_{t_{i-1}^n} ; \theta)$, which is
a martingale estimating function. The
close relationship between martingale estimating functions and 
the type of GMM-estimators described here is discussed in detail in
\cite{BJ}. More general GMM-estimators of the martingale estimating
function type were considered in \cite{hansen2, hansen93}, and
a discussion of links between the literature on estimating functions 
and that on GMM-estimators can be found in \cite{hansen2001}.

Approximate martingale estimating functions can be obtained by 
replacing the exact conditional expectation in (\ref{specielestfct}) 
by the approximation given by (\ref{expansion}) such that the function
$g$ has the form
\be
\label{specielestfctapprox}
g(\Delta,y,x; \theta) = A(x, \Delta;\theta) \left[f(y;\theta) -
\pi_\theta^{\kappa,\Delta} f (x;\theta) \right], 
\ee
where
\be
\label{transitionopapprox}
\pi_\theta^{\kappa,\Delta} f (x;\theta) = 
\sum_{i=0}^{\kappa-1} \frac{\Delta^i}{i!} \, L_\theta^i f (x;\theta),
\hspace{6mm} \kappa = 2,3,\ldots,
\ee
with the generator $L_\theta$ applied coordinate-wise.
This estimating functions satisfies (\ref{mart}). 
A simple example is $g(\Delta,y,x;\theta) = a(x,\Delta;\theta)
(y - x - b(x;\alpha)\Delta)$ with $\kappa = 2$, considered by
\cite{rao2} and \cite{zmirou}. Other instances are
the estimators proposed by \cite{ckls} and \cite{kelly}.
For all $\kappa \in \N$, ($\kappa \geq 2$), \cite{kessler97} proposed
a {\it Gaussian approximation} to the likelihood function, for 
which the corresponding pseudo-score function is an approximate 
martingale estimating function that satisfies (\ref{mart}).

\sect{Optimal rate}
\label{ratesection}

In this section we present asymptotic results for 
approximate martingale estimating functions. We begin with a general
approximate martingale estimating function. Then we will see how
Condition \ref{con3} implies rate optimality, so 
that the estimator of the parameter in the diffusion 
coefficient converges faster than the estimator of the parameter in the 
drift coefficient. As previously, $x^T$ denotes the 
transpose of a vector or matrix $x$ and $\theta_0 =
(\alpha_0,\beta_0)$ denotes the true parameter value. 

\begin{theo}
\label{theorem1}
Assume that the Conditions \ref{con1} and \ref{con2} hold. Suppose,
moreover, the identifiability condition that
\bea
\label{gamma}
\lefteqn{\gamma(\theta,\theta_0) =
\int_\ell^r [b(x,\alpha_0)-b(x,\alpha)]\partial_y g(0,x,x;\theta)
\mu_{\theta_0} (x) dx} \\
&& \hspace{18mm}\mbox{} + \mbox{\small $\frac12$} \int_\ell^r 
[v(x,\beta_0)-v(x,\beta)]\partial^2_y g(0,x,x;\theta)
\mu_{\theta_0} (x) dx \neq 0 \nonumber
\eea
for all $\theta \neq \theta_0$, and that the matrix
\be
\label{S}
S = \int_\ell^r J_{\theta_0}(x) \mu_{\theta_0}(x) dx
\ee
is invertible, where
\be
\label{A}
J_\theta(x) = \left( \begin{array}{cc}
\partial_\alpha b(x;\alpha)
\partial_y g_1(0,x,x ; \theta) &
\mbox{\small $\frac12$} \partial_\beta v(x;\beta)
\partial^2_y g_1(0,x,x ; \theta) \\ & \\
\partial_\alpha b(x;\alpha)
\partial_y g_2(0,x,x ; \theta) &
\mbox{\small $\frac12$} \partial_\beta 
v(x;\beta)
\partial^2_y g_2(0,x,x ; \theta)
\end{array} \right).
\ee
Then with a probability that goes to 
one as $n \rightarrow \infty$, a consistent $G_n$-estimator $\hat{\theta}_n
= (\hat{\alpha}_n,\hat{\beta}_n)$ exists and is unique in any compact subset 
$K$ of $\Theta$ with $\theta_0 \in$ int$\, K$. If $n \Delta_n^{2\kappa-1} \rightarrow
0$, then
\be
\label{asnorm1}
\sqrt{n \Delta_n} (\hat{\theta}_n - \theta_0) 
\stackrel{\cal D}{\longrightarrow} 
N_2 \left( 0, S^{-1} V_0 (S^T)^{-1} \right)
\ee
under $P_{\theta_0}$, where $V_0 = V(\theta_0)$ with
\be
\label{Vtheta}
V(\theta) = \int_\ell^r v(x,\beta_0) \partial_y g(0,x,x;\theta)
\partial_y g(0,x,x;\theta)^T \mu_{\theta_0} (x) dx.
\ee 
For a martingale estimating function (\ref{asnorm1}) holds without the
extra condition on the rate of convergence of $\Delta_n$.

A consistent estimator of the asymptotic variance of $\hat{\theta}_n$
can be obtained from
\be
\label{varest1}
\frac{1}{n\Delta_n} \sum_{i=1}^n \partial_{\theta^T} 
g(\Delta_n, X_{t_i^n}, X_{t_{i-1}^n}; \hat{\theta}_n) 
\stackrel{P_{\theta_0}}{\longrightarrow} -S
\ee
\be
\label{varest2}
\frac{1}{n\Delta_n} \sum_{i=1}^n g(\Delta_n, X_{t_i^n}, X_{t_{i-1}^n};
\hat{\theta}_n) g(\Delta_n, X_{t_i^n}, X_{t_{i-1}^n}; \hat{\theta}_n)^T 
\stackrel{P_{\theta_0}}{\longrightarrow} V_0
\ee

\end{theo}

The theorem follows from general asymptotic statistical results for stochastic
processes, see e.g.\ \cite{jjms}. The proof is given in Section
\ref{proofs}. The precise meaning of the uniqueness statement is that
for any $G_n$-estimator $\tilde \theta_n$ with $P_{\theta_0}(\tilde \theta_n
\in K) \rightarrow 1$ as $n \rightarrow \infty$, it holds that
$P_{\theta_0}(\hat \theta_n \neq \tilde \theta_n) \rightarrow 0$ as $n
\rightarrow \infty$. 

We see from (\ref{asnorm1}) that the rate of convergence of both 
$\hat{\alpha}$ and $\hat{\beta}$ is $1/\sqrt{n \Delta_n}$, provided
that the matrix $V_0$ is regular. Here $n \Delta_n$ is the length of
the interval in which the diffusion is observed. 
\cite{gobet} showed that under weak regularity conditions a discretely
sampled diffusion model is local asymptotically normal in the high
frequency/infinite time horizon  asymptotic scenario considered here,
and that the optimal rate of convergence for estimators of parameters in
the drift coefficient is $1/\sqrt{n \Delta_n}$, whereas the 
optimal rate for estimators of parameters in the diffusion coefficient
is $1/\sqrt{n}$. 

The next theorem shows what happens when Jacobsen's condition, 
Condition \ref{con3}, is satisfied, or more precisely, when a version of the
estimating function satisfies the condition. 

\begin{theo}
\label{theorem2}
Suppose the Conditions \ref{con3}, \ref{con1} and \ref{con2} hold. 
Assume, moreover, that the following identifiability condition is 
satisfied
\bean
\int_\ell^r [b(x,\alpha_0)-b(x,\alpha)]\partial_y g_1(0,x,x;\theta)
\mu_{\theta_0} (x) dx & \neq & 0 \hspace{6mm} \mbox{when } 
\alpha \neq \alpha_0 \\ \\
\int_\ell^r [v(x,\beta_0)-v(x,\beta)]\partial^2_y g_2(0,x,x;\theta)
\mu_{\theta_0} (x) dx & \neq & 0 \hspace{6mm} \mbox{when } 
\beta \neq \beta_0,
\eean
and that $S_{11} \neq 0$ and $S_{22} \neq 0$ with $S$ given by
(\ref{S}). 
Then with a probability that goes to one as $n \rightarrow
\infty$, a consistent $G_n$-estimator $\hat{\theta}_n = (\hat{\alpha}_n,
\hat{\beta}_n)$ exists and is  unique in any compact subset K of
$\Theta$ with $\theta_0 \in$  int$\, K$. 

If, moreover,
\be
\label{extracond}
\partial_\alpha \partial^2_y g_2 (0,x,x;\theta) = 0
\ee
and $n \Delta_n^{2(\kappa-1)} \rightarrow 0$, then
\be
\label{asnorm2}
\left( \begin{array}{c} \sqrt{n \Delta_n} (\hat{\alpha}_n - \alpha_0)
\\  \sqrt{n} (\hat{\beta}_n - \beta_0) \rule{0mm}{6mm} 
\end{array} \right) \stackrel{\cal D}{\longrightarrow} 
N_2 \left( \left( \begin{array}{c} 0 \\ 0 \rule{0mm}{4mm} \end{array} \right), 
\left( \begin{array}{cc} W_1(\theta_0)/S^2_{11} & 0 \\  0 & 
W_2(\theta_0)/S^2_{22} \rule{0mm}{6mm} \end{array} \right)  \right)
\ee
where
\bean
W_1(\theta) &=& \int_\ell^r v(x;\beta_0)[\partial_y g_1 
(0,x,x;\theta)]^2 \mu_{\theta_0} (x) dx = V(\theta)_{11} \\
W_2(\theta) &=& \mbox{\small $\frac12$} \int_\ell^r [v(x; \beta_0)^2
+  \frac12 (v(x; \beta_0) - v(x; \beta))^2]
[\partial^2_y g_2 (0,x,x;\theta)]^2 \mu_{\theta_0}(x) dx
\eean
with $V(\theta)$ given by (\ref{Vtheta}). For a martingale estimating
function (\ref{asnorm2}) holds without the extra condition on the rate
of convergence of $\Delta_n$. 

A consistent estimator of the asymptotic variance of $\hat{\theta}_n$
can be obtained from (\ref{varest1}) and
\be
\label{varest3}
D_n \sum_{i=1}^n g(\Delta_n, X_{t_i^n}, X_{t_{i-1}^n};
\hat{\theta}_n) g(\Delta_n, X_{t_i^n}, X_{t_{i-1}^n}; \hat{\theta}_n)^T D_n 
\stackrel{P_{\theta_0}}{\longrightarrow}
\left( \begin{array}{cc} W_1(\theta_0) & 0 \\ 0 & W_2(\theta_0) \rule{0mm}{6mm}
\end{array} \right),
\ee
where
\be
\label{Dn}
D_n = \left( \begin{array}{cc} \frac{1}{\sqrt{n\Delta_n}} & 0 \\ 
0 & \frac{1}{\Delta_n \sqrt{n}} \rule{0mm}{6mm} \end{array} \right).
\ee
\end{theo}

\vspace{4mm}

Thus Jacobsen's condition (\ref{jacobsen}) and the additional
condition (\ref{extracond}) imply rate optimal estimators and that the
estimators of the drift parameter and of the diffusion coefficient
parameter are asymptotically independent. In the next section 
we shall see that (\ref{extracond}) is automatically satisfied for
efficient estimating functions. Note that for non-martingale
estimating functions $\Delta_n$ must go a bit faster to zero than was
required in Theorem \ref{theorem1}. Note also that if the
first coordinate of $g$ satisfies Jacobsen's condition too, then the
first part of the identifiability condition in Theorem \ref{theorem2} 
does not hold, and the
parameter $\alpha$ cannot be consistently estimated by the estimating
function (\ref{estfct}). The proof of Theorem \ref{theorem2} is given
in Section \ref{proofs}.

\vspace{3mm}

\begin{ex}
\label{quadratic0}
{\em Consider a} quadratic martingale estimating function 
{\em of the form
\be
\label{quad}
g(\Delta,y,x;\theta) =  \left( \begin{array}{c} a_1(x,\Delta;\theta)
[y-F(\Delta,x;\theta)] \\
a_2(x, \Delta; \theta) \left[ (y-F(\Delta,x; \theta))^2 - 
\phi(\Delta,x;\theta) \right] \rule{0mm}{6mm} \end{array} \right),
\ee
where $F(\Delta,x;\theta) = E_{\theta} (X_{\Delta}| X_0 = x)$ and
$\phi(\Delta,x;\theta) = \mbox{Var}_{\theta} (X_{\Delta}| X_0 = x)$.
Since, by (\ref{expansion}), $F(\Delta,x;\theta) =
x + O(\Delta)$ and $\phi(\Delta,x;\theta) = O(\Delta)$, we find that
\be
\label{quad0}
g(0,y,x;\theta) = \left( \begin{array}{c}
a_1(x,0;\theta)(y-x) \\
a_2(x,0;\theta)(y-x)^2 \rule{0mm}{6mm}
\end{array} \right).  
\ee
Jacobsen's condition (\ref{jacobsen}) is satisfied because $\partial_y
g_2(0,y,x;\theta) = 2 a_2(x,\Delta;\theta)(y-x)$. Thus estimators 
obtained from (\ref{quad}) are rate optimal, provided that 
(\ref{extracond}) is satisfied, for instance if $a_2$ does 
not depend on $\alpha$.

Clearly (\ref{quad0}) holds if $F$ and $\phi$ in (\ref{quad}) are
replaced by expansions of order $O(\Delta^{\kappa-1})$ with $\kappa
\geq 2$, using again  (\ref{expansion}). Thus rate optimal estimators
are also obtained in this more easily calculated case, provided again that 
(\ref{extracond}) holds. The simplest example ($\kappa=2$) is
\be
\label{euler}
g(\Delta,y,x;\theta) =  \left( \begin{array}{c} a_1(x,\Delta;\theta)
[y-x-b(x;\alpha)\Delta] \\
a_2(x, \Delta; \theta) \left[ (y-x-b(x; \alpha)\Delta)^2 - 
v(\Delta,x;\beta)\Delta \right] \rule{0mm}{6mm} \end{array} \right).
\ee

It is instructive to consider an example of an estimating
function for which estimators are not rate optimal. The martingale
estimating function 
\be
\label{non-rate}
g(\Delta,y,x;\theta) =  \left( \begin{array}{c} a_1(x,\Delta;\theta)
[y-F(\Delta,x;\theta)] \\
a_2(x, \Delta; \theta) \left[ y^2-(\phi(\Delta,x;\theta) + 
F(\Delta,x;\theta)^2) \right] \rule{0mm}{6mm} \end{array} \right)
\ee
does not satisfy (\ref{jacobsen}). It is easy to check that a version
of $g$ satisfying (\ref{jacobsen}) exists if and only if
$a_1(x,0;\theta) = c_\theta a_2(x,0;\theta)x$ for some real constant
$c_\theta$. In all other cases, the estimating function given by (\ref{non-rate}) is 
not rate optimal. We can obtain a particular case of (\ref{quad}) from
(\ref{non-rate}) by choosing $a_1(x,\Delta;\theta) =
a_2(x,\Delta;\theta)F(\Delta,x;\theta)$ and using the version $\tilde
g_1 = g_1, \ \tilde g_2 = g_2 - 2g_1$. This particular case of
(\ref{non-rate}), obviously satisfies that $a_1(x,0;\theta) = a_2(x,0;\theta)x$.
}

\qed

\end{ex}

\sect{Efficient estimating functions}
\label{efficientsection}

In this section we study the conditions under which an approximate
martingale estimating function, $G_n(\theta)$, yields an efficient
estimator. In particular, we show that Condition \ref{conefficiency}  
ensures efficiency, and that Godambe-Heyde optimal estimating
functions yield rate optimal and efficient estimators.

\begin{theo}
\label{theorem:efficiency}
Suppose Condition \ref{conefficiency} and the conditions of Theorem
\ref{theorem2} except (\ref{extracond}) are satisfied. 
Then the conclusions of Theorem \ref{theorem2} hold and the estimating
function (\ref{estfct}) is efficient, i.e.,  the asymptotic covariance 
matrix of the estimator $\hat{\theta}_n = (\hat{\alpha}_n,\hat{\beta}_n)$ 
equals
\be
\label{optimalcovariance}
\Sigma(\theta_0) =
\left( \begin{array}{cc} \left( {\displaystyle \int_\ell^r} 
\frac{(\partial_\alpha  b(x;\alpha_0))^2}{v(x;\beta_0)}\mu_{\theta_0}(x)dx 
\right)^{-1} & 0 \\
0 & 2 \left( {\displaystyle \int_\ell^r} 
\left[\frac{\partial_\beta v(x;\beta_0)}{v(x;\beta_0)}\right]^2
\mu_{\theta_0}(x)dx \right)^{-1} \rule{0mm}{6mm} \end{array} \right).
\ee
Consistent estimators of the asymptotic variances can be obtained from
\[
\frac{1}{n \Delta_n} \sum_{i=1}^n g_1(\Delta_n, X_{t_i^n}, X_{t_{i-1}^n};
\hat \theta_n )^2 \stackrel{P_{\theta_0}}{\longrightarrow} \int_\ell^r 
\frac{(\partial_\alpha  b(x;\alpha_0))^2}{v(x;\beta_0)}\mu_{\theta_0}(x)dx
\]
and
\[
\frac{1}{n \Delta_n^2} 
\sum_{i=1}^n g_2(\Delta_n, X_{t_i^n}, X_{t_{i-1}^n}; \hat \theta_n)^2 
\stackrel{P_{\theta_0}}{\longrightarrow} \int_\ell^r
\left[\frac{\partial_\beta v(x;\beta_0)}{v(x;\beta_0)}\right]^2
\mu_{\theta_0}(x)dx.
\]
\end{theo}

Note that an efficient estimating function automatically
satisfies (\ref{extracond}). An asymptotic martingale estimating
function is efficient if and only if there exists a version that
satisfies the conditions of Theorem \ref{theorem:efficiency}.

The covariance matrix (\ref{optimalcovariance}) is equal to the
leading term in the expansion of the asymptotic 
variance of the  maximum likelihood estimator in powers of $\Delta$
found by \cite{castellezmirou}. The asymptotic variance of $\hat \alpha_n$ 
equals that of the maximum likelihood estimator based on continuous time
observation, see e.g.\ \cite{kutoyants}. 

\begin{ex}
\label{quadratic}
{\em Consider again the} quadratic martingale estimating function 
{\em (\ref{quad}). The function $g(0,y,x;\theta)$, 
given by (\ref{quad0}), satisfies the conditions
for efficiency (\ref{drifteff}) and (\ref{diffeff}) if we choose
$a_1(x,\Delta;\theta) = \partial_\alpha  b(x; \alpha)/\sigma^2(x; \beta)$
and $a_2(x,\Delta;\theta) = \partial_\beta \sigma^2(x; \beta)/\sigma^4(x; \beta)$,
as proposed by \cite{bmb&ms1, bibbyso96}. The same is true
of weight functions $a_1$ and $a_2$ that converge to
$\partial_\alpha b/\sigma^2$ and $\partial_\beta \sigma^2/\sigma^4$ as $\Delta \rightarrow 0$.
An example is $a_1(x,\Delta;\theta) = \partial_\alpha
F(\Delta,x;\theta)/\phi(\Delta,x;\theta)$  
and $a_2(x,\Delta;\theta) = \Delta \partial_\beta \phi(\Delta, x;\theta)/\phi(\Delta,
x;\theta)^2$. This is the optimal quadratic martingale estimating
function in the sense of \cite{godambeheyde} (after multiplication of
the second coordinate by $\Delta$), see \cite{bmb&ms1, bibbyso96}. 

Consider the pseudo-likelihood function obtained from the likelihood
function by replacing the transition density $p(\Delta,y,x;\theta)$ by the
Gaussian density with mean $F(\Delta,x;\theta)$
and variance $\phi(\Delta,x;\theta)$. The exact conditional moments
are used to ensure consistency of the estimator also in the
case of low frequency asymptotics, where $\Delta$ is not small. The
corresponding pseudo-score function is the quadratic estimating function 
(\ref{quad}) with $a_1(x,\Delta;\theta) = \partial_\alpha
F(\Delta,x;\theta)/\phi(\Delta,x;\theta)$  
and $a_2(x,\Delta;\theta) = \partial_\beta \phi(\Delta, x;\theta)/\phi(\Delta,
x;\theta)^2$, which we have just seen is efficient.

A pseudo-likelihood function that works for data sampled at a high
frequency is the likelihood function obtained by replacing the
original diffusion model by its Euler approximation. It can be
obtained from the original likelihood function by replacing the
transition density by the Gaussian density with mean and variance
given by the expansions $x-b(x;\alpha)\Delta$ and $\sigma^2(x;\beta)\Delta$.
The corresponding pseudo score is of the form (\ref{euler}) with 
$a_1(x,\Delta;\theta) = \partial_\alpha  b(x; \alpha)/\sigma^2(x; \beta)$
and (after multiplication by $\Delta$) $a_2(x,\Delta;\theta) =
\partial_\beta \sigma^2(x; \beta)/\sigma^4(x; \beta)$. Since  (\ref{quad0}) holds,
this Euler pseudo score function satisfies the conditions for
efficiency.
This estimator has often been used in empirical work in finance.
Similarly, it follows that the estimators considered by 
\cite{dorogovcev}, \cite{rao2}, \cite{zmirou}, \cite{nakahiro3}, 
\cite{kessler97}, \cite{kelly}, and \cite{uchidayoshida} are efficient
under suitable conditions on the rate of convergence of $\Delta_n$.
}

\qed

\end{ex}

\begin{ex}
\label{mleex}
{\em A final example is} maximum likelihood estimation. {\em In broad
generality, the score function is a martingale estimating
function, see e.g.\ \cite{oebnms}.
The transition density can, under weak regularity conditions, be
expanded in powers of $\Delta$
\[
p(\Delta,y,x;\theta) = r(\Delta,y,x;\theta)(1 + O(\Delta)),
\]
where
\bean
\lefteqn{r(\Delta,y,x;\theta) =} \\ &&  \frac{1}{\sqrt{2 \pi \sigma^2(y;\beta)\Delta}} 
\exp \left( -\frac{(k(y;\beta)-k(x;\beta))^2}{2 \Delta}
  +m(y; \alpha,\beta) - m(x; \alpha,\beta)  -\mbox{\small $\frac12$}
  \log \left(\frac{\sigma (y;\beta)}{\sigma (x;\beta)} \right) \right),
\eean
$k(x;\beta) = \int^x \sigma^{-1}(z;\beta)dz$ and $m(x; \alpha,\beta) =
\int^x b(z;\alpha)/\sigma^2(z;\beta)dz $, 
see e.g.\ \cite{castellezmirou} or \cite{gihman}, Part I, Chapter
13. Therefore, under regularity conditions, the (suitably normalized) score function 
$g_1(\Delta,y,x;\theta) = \partial_\alpha \log p(\Delta,y,x;\theta)$
and $g_2(\Delta,y,x;\theta) = \Delta \partial_\beta \log p(\Delta,y,x;\theta)$
satisfies that
\bean
g_1(\Delta,y,x;\theta) &=&  \int_x^y \frac{\partial_\alpha b(z;\alpha)}
{\sigma^2(z;\beta)}dz + O(\Delta) \\ && \\
g_2(\Delta,y,x;\theta) &=&   -[k(y;\beta) - k(x;\beta)]
[\partial_\beta k(y;\beta) - \partial_\beta k(x;\beta)] + O(\Delta).
\eean
From these expansions it follows that the score functions
(normalized as above) satisfies the conditions (\ref{jacobsen}),
(\ref{drifteff}) and (\ref{diffeff}) for rate optimality and efficiency.
In particular, $\partial^2_y g_2(0,x,x;\theta) = - 2 \partial_x
k(x;\beta) \partial_\beta \partial_x k(x;\beta) = 
\partial_\beta \sigma^2(x;\beta)/\sigma^4(x;\beta)$. Obviously,
also the pseudo-likelihood function obtained by replacing the
transition density $p$ by $r$ has a pseudo-score
function that satisfies the conditions for rate optimality and efficiency.}

\qed

\end{ex}

Consider martingale estimating functions of the form 
(\ref{specielestfct}) and the related approximate martingale estimating 
functions (\ref{specielestfctapprox}), i.e.
\be
\label{compact}
G_n(\theta) = \sum_{i=1}^n A(X_{t_{i-1}^n}, \Delta; \theta)
[f(X_{t_{i}^n};\theta) - \pi_\theta^{\kappa,\Delta} 
f( X_{t_{i-1}^n};\theta)], \hspace{5mm} \kappa = 1,2,\ldots,
\ee
with $\pi_\theta^{\kappa,\Delta}$ given by (\ref{transitionop}) for
$\kappa = 1$ (the martingale case) and by (\ref{transitionopapprox})
for $\kappa = 2,3,\ldots$. Moreover, $A$ is a $2 \times N$-matrix of
weights, and we assume that the coordinates of the $N$-dimensional
function $f(x,\theta)$ are twice continuously differentiable w.r.t. $x$. 

For estimating functions of the form (\ref{compact}),
the condition for rate optimality is  
\be
\label{rate44}
\sum_{j=1}^N a_{2j}(x,0;\theta)\partial_x f_j(x;\theta) = 0, 
\ee
where $a_{ij}$ denotes the $ij$th entry of $A$, and the condition 
for efficiency is
\bea
\sum_{j=1}^N a_{1j}(x,0;\theta)\partial_x f_j(x;\theta) &=& 
\partial_\alpha b(x;\alpha)/\sigma^2(x;\beta) \label{eff44a}\\
\sum_{j=1}^N a_{2j}(x,0;\theta)\partial^2_x f_j(x;\theta)
&=& \partial_\beta \sigma^2(x;\beta)/\sigma^4(x;\beta)
\label{eff44b}
\eea
For a given function $f$, we want to find a weight-matrix $A$ such
that these equations are satisfied. Obviously, it is necessary that
$N \geq 2$ in order that all three equations are satisfies.  If $N=1$,
an efficient approximate martingale estimating function can be
obtained by solving (\ref{eff44a}), provided that the diffusion
coefficient is known, so that only the drift depends on a parameter.

First consider $N = 2$, and assume that the matrix
\be
\label{D}
M(x) = \left( \begin{array}{cc}
\partial_x f_1(x;\theta) & \partial^2_x f_1(x;\theta) \\ 
\partial_x f_2(x;\theta) & \partial^2_x f_2(x;\theta) \rule{0mm}{6mm} 
\end{array} \right)  
\ee
is invertible for $\mu_\theta$-almost all $x$. Then the linear
equations (\ref{rate44}) - (\ref{eff44b}) are satisfied for
\be
\label{deltaoptimal}
A(x, 0 ;\theta) = \left( \begin{array}{cc}
\partial_\alpha b(x;\alpha)/\sigma^2(x;\beta) & c(x ; \theta) \\
0 & \partial_\beta \sigma^2(x;\beta)/\sigma^4(x;\beta) \rule{0mm}{7mm}
\end{array} \right) M(x)^{-1},
\ee
where $c(x ; \theta)$ is any (measurable) function. As a simple example,
the quadratic estimating function ($N=2$, $f_1(x)=x$ and
$f_2(x)=\frac12 x^2$) is rate optimal and efficient if the weights are
weights $a_{11}(x) = \partial_\alpha b(x; \alpha)/\sigma^2(x;\beta)$,
$a_{12}(x) = 0$, $a_{22}(x) = \partial_\beta
\sigma^2(x;\beta)/\sigma^4(x;\beta)$ and $a_{21}(x) = -2xa_{22}(x)$ 
(we have chosen $c = 0$). Note that a simple choice for
the weight matrix $A(x,\Delta;\theta)$ for $\Delta > 0$ is to let it
be given by (\ref{deltaoptimal}) for all $\Delta$.

It is easily seen that, for any $N \geq 2$, there exist many solutions to
(\ref{rate44}) - (\ref{eff44b}) provided that there are two coordinates
of $f$ (without loss of generality,  we can assume these to be $f_1$
and $f_2$) such that  $M(x)$ is invertible. In the special case
$\kappa =1$, this result follows from Theorem 2.2 of \cite{mj1}.
The conditions for Jacobsen's concept small $\Delta$-optimality for
martingale estimating functions are identical to our conditions for
rate optimality and efficiency, so we can take advantage of his thorough
study of when martingale estimating functions of the type
(\ref{compact}) with $\kappa = 1$ are small $\Delta$-optimal.  

Here we will, however, go another way and give a natural and generally
useful way of finding a rate optimal and efficient weight matrix $A$
for all $\kappa \geq 1$. A weight matrix $A^*$ in a martingale
estimating function of the type (\ref{compact}) (i.e.\ $\kappa =1$) is
optimal in the sense of \cite{godambeheyde}, see also \cite{heyde97},
if it solves the linear equation 
\bea
\label{godambeheyde}
A^*(x,\Delta;\theta) & \hspace{-15mm} E_\theta \left( [f(X_\Delta;\theta) -
\pi_\theta^{1,\Delta} f(x;\theta)] [f(X_\Delta;\theta) -
\pi_\theta^{1,\Delta} f(x;\theta)]^T \, | \, X_0 = x \right) \\
& \hspace{60mm} \mbox{}
= \partial_\theta  \pi_\theta^{1,\Delta} f^T(x;\theta) - 
\pi_\theta^{1,\Delta} \partial_\theta f^T(x;\theta). \nonumber
\eea
It can be assumed that the functions $f_1, \ldots, f_N$ are affinely
independent such that the conditional covariance matrix in
(\ref{godambeheyde}) is invertible.  
The Godambe-Heyde optimal martingale estimating function gives
an estimator that minimizes the asymptotic variance of estimators
obtained from the class of martingale estimating functions of the form
(\ref{compact}) (with $\kappa = 1$) with a fixed function $f$ and
for a fixed, possibly large, $\Delta$. The next theorem shows that the
Godambe-Heyde optimal estimators are rate optimal and 
efficient in the high frequency asymptotic scenario considered in the
present paper. Moreover, the same is true of the approximate
martingale estimating functions obtained by expanding all conditional
moments (including those in $A^*$) in powers of
$\Delta$, which gives a feasible general way of constructing
explicit estimating functions that are rate optimal and efficient.

\begin{theo}
\label{theorem4}
Suppose Condition \ref{con1} is satisfied, that $f_j \in 
C_{p,6,1}((\ell,r)\times \Theta)$, 
$j=1,\ldots,N$, that $N \geq 2$ and
that the $2 \times 2$ matrix $M(x)$ given by (\ref{D}) is invertible
for $\mu_\theta$-almost all $x$. Let $A^* (x, \Delta; \theta)$ satisfy
(\ref{godambeheyde}), and define 
\[
B(x,\Delta,,\theta)= \left( \begin{array}{cc}
1 & 0 \\ 0 & 2 \Delta \rule{0mm}{5mm} \end{array} \right)  
A^* (x, \Delta; \theta).
\]
Then the limit $B(x,0,\theta)$ exists, and
\be
\label{gstar}
g^*(\Delta,y,x;\theta) = B(x,\Delta,,\theta)  
[f(y;\theta) - \pi_\theta^{\kappa,\Delta} f(x;\theta)]
\ee
satisfies the conditions for rate optimality (\ref{jacobsen}) 
and efficiency (\ref{drifteff}) and (\ref{diffeff}) for all 
$\kappa \in \N$. The same is true
if $B$ is replaced by a matrix $\tilde B$ satisfying that $\tilde 
B(x,0,\theta) = B(x,0,\theta)$.
\end{theo}

The matrix $\tilde B$ can be obtained by replacing the conditional
moments in $A^*$ by expansions in powers of $\Delta$, or simply by
defining $\tilde B(x,\Delta,\theta) = B(x,0,\theta)$. It is surprising
that a local property like Godambe-Heyde optimality,  
which ensures optimality only within a particular class of estimating
functions, implies global optimality properties like rate optimality and efficiency.
Phrased in terms of the concept small $\Delta$-optimality, this
result was conjectured by \cite{mj1} for martingale estimating  
functions ($\kappa =1$).
The fact that one of the conditions for efficiency is $N \geq 2$ explains 
the finding in \cite{larsen}
that an optimal martingale estimating function based on two 
eigenfunctions seemed to be efficient for weekly observations of
exchange rates in a target zone.

Let us conclude this section by stating the results for a $d$-dimensional
diffusion. In this case $b(x;\alpha)$ is $d$-dimensional and $v(x;\beta)
= \sigma(x;\beta)\sigma(x;\beta)^T$ is a $d \times d$-matrix. The
conditions for efficiency are
\[
\partial_y g_1 (0,x,x;\theta) = \partial_\alpha b(x;\alpha)^T v(x;\beta)^{-1}  
\]
and
\[
\mbox{vec}\left(\partial^2_y g_2 (0,x,x;\theta)\right) 
= \mbox{vec}\left(\partial_\beta 
v(x;\beta)\right) \left(v^{\otimes 2}(x;\beta)\right)^ {-1}.
\]
In the latter equation, vec$ (M)$ denotes for a $d \times d$ 
matrix $M$ the $d^2$-dimensional row vector consisting of the 
rows of $M$ placed one after the other, and $M^{\otimes 2}$ is the
$d^2 \times d^2$-matrix with $(i',j'),(ij)$th entry equal to
$M_{i'i}M_{j'j}$. Thus if $M = \partial_\beta 
v(x;\beta)$ and $M^\bullet = (v^{\otimes 2}(x;\beta))^ {-1}$, then
the $(i,j)$th coordinate of vec$ (M) \, M^\bullet$ is $\sum_{i'j'}
M_{i'j'} M^\bullet_{(i'j'),(i,j)}$. These expressions are the 
conditions for small
$\Delta$-optimality for multivariate diffusions given by \cite{mj1}.

For a $d$-dimensional diffusion process, the condition analogous to
the one discussed shortly before Theorem \ref{theorem4} ensuring the
existence of a rate optimal and efficient estimating function of
the form (\ref{compact}) is that $N \geq d(d+3)/2$, and
that the $N \times (d + d^2)$-matrix
\[
\left( \begin{array}{cc}
\partial_x f(x;\theta) &  \partial^2_x f(x;\theta) \end{array} \right) 
\] 
has full rank $d(d+3)/2$. For $\kappa =1$ this follows from Theorem
2.2 of \cite{mj1}, and it is clear from the proof of this theorem that
it holds for $\kappa \geq 2$ too. When $\alpha$ and $\beta$
are multivariate, we further need that $\{ \partial_{\alpha_i}
b(x;\alpha) \}$ and $\{ \partial_{\beta_i} v(x;\beta) \}$ are two sets of
linearly independent functions of $x$. These conditions also ensure
that Theorem \ref{theorem4} holds for a  $d$-dimensional
diffusion process, i.e.\ that the Godambe-Heyde optimal martingale
estimating function is rate optimal and efficient for a  $d$-dimensional
diffusion process.

\sect{Proofs and lemmas}
\label{proofs}

The first of the following lemmas is a slight generalization of Lemma 6 in
\cite{kessler97}, while the second lemma is essentially Lemma 8 in the same
paper. The proofs are analogous to those in Kessler's paper. The result 
(\ref{lem5.3}) follows from (\ref{lem5.2}). The notation
$R(\Delta,y,x;\theta)$ was defined in Section \ref{model}. We sometimes use
the notation $a \leq_C b$, which means that there exists a $C>0$ such
that $a \leq Cb$.

\begin{lem}
\label{lemma5}
Assume Condition \ref{con1}. Then a constant $C_k > 0$ exists for $k =
1,2, \ldots$ such that 
\be
\label{lem5.1}
E_{\theta_0}( |X_{t+\Delta}-X_t|^k \, | \, X_t) \leq C_k \Delta^{k/2}
(1+|X_t|)^{C_k}
\ee
for $\Delta > 0$. Let $f(y,x,\theta)$ be a real function of 
polynomial growth in $x$ and $y$ uniformly for $\theta$ in a compact 
set $K$. Then for any fixed $\Delta_0 >0$ there exists a constant $C >
0$ such that 

\be 
\label{lem5.2}
E_{\theta_0}( |f(X_{t+\Delta},X_t, \theta)| \, | \, X_t) 
\leq C (1+|X_t|)^C \ \ \ \mbox{ for } \Delta \in [0,\Delta_0] 
\mbox{ and } \theta  \in K.
\ee
Suppose the function $f(y,x,\theta)$
is, moreover, $2k$ times differentiable ($k \leq 3$) 
with respect to $y$ with all derivatives of polynomial growth in $x$ 
and $y$ uniformly for $\theta$ in compact sets. Then
\be
\label{lem5.3}
\int_0^\Delta  \int_0^{u_1} \cdots 
\int_0^{u_{k-1}} E_{\theta_0}\left( L_{\theta_0}^{k} (f)
  (X_{t+u_{k}},X_t; \theta) \, | \, X_t \right) du_{k} \cdots du_1 
= \Delta^{k} R(\Delta, X_t,\theta).
\ee
\end{lem}

\vspace{4mm}

The result (\ref{lem5.3}) is used to ensure that the remainder term 
in expansions of the type (\ref{expansion}) have the expected order. 
The result can be proved for larger values of $k$ if stronger differentiability
conditions are imposed on the coefficients $b$ and $\sigma$.

\vspace{3mm}

\begin{lem}
\label{lemma4}
Assume Condition \ref{con1}, and
let $f(x,\theta)$ be a real function that is differentiable with respect to
$x$ and $\theta$ with derivatives of polynomial growth in $x$ 
uniformly for $\theta$ in a compact set. Then
\[
\frac1n \sum_{i=1}^n f(X_{t_{i}^n};\theta) 
\stackrel{P_{\theta_0}}{\longrightarrow} 
\int_\ell^r f(x;\theta) \mu_{\theta_0} (x) dx
\]
uniformly for $\theta$ in  a compact set.
\end{lem}

\vspace{3mm}

Lemma 9 in \cite{gcjj} is used frequently in the proofs of 
Lemma \ref{lemma2} and Lemma \ref{lemma3} to establish
pointwise convergence. The result is therefore cited here for
the convenience of the reader.

\vspace{3mm}

\begin{lem}
\label{lemmagcjj} Let $Z^n_i$ ($i=1, \ldots, n, n \in \N$) be a triangular 
array of random variables such that $Z^n_i$ is 
${\cal G}^n_{i}$-measurable, where ${\cal G}^n_{i} = \sigma (W_s \, 
: \, s \leq t^n_i)$. If
\[
\sum_{i=1}^n E_\theta(Z^n_i \, | \, {\cal G}^n_{i-1}) 
\stackrel{P_{\theta}}{\longrightarrow} U
\]
and
\[
\sum_{i=1}^n E_\theta((Z^n_i)^2 \, | \, {\cal G}^n_{i-1}) 
\stackrel{P_{\theta}}{\longrightarrow} 0,
\]
where $U$ is a random variable, then
\[
\sum_{i=1}^n Z^n_i \stackrel{P_{\theta}}{\longrightarrow} U.
\]

\end{lem}

\noindent
{\it Proof of Lemma \ref{lemma1}.} Combining (\ref{gdeltaexp})
and (\ref{expansion}), we find that
\[
E_\theta (g(\Delta_n,X_{t_{i}^n},X_{t_{i-1}^n};
\theta) \, | \, X_{t_{i-1}^n}) \\
= \sum_{\ell =0}^{\kappa-1} \frac{\Delta_n^\ell}{\ell !}
\sum_{j=0}^\ell {\ell \choose j} L_\theta^{\ell -
  j}(g^{(j)}(\theta))(X_{t_{i-1}^n},X_{t_{i-1}^n}) 
+ \Delta_n^{\kappa} R(\Delta, X_{t_{i-1}^n},\theta),
\]
from which the ``if'' statement of the lemma follows immediately. The
``only if'' statement follows from the same expansion because an
approximate martingale estimating function satisfies $E_\theta
(g(\Delta_n,X_{t_{i}^n},X_{t_{i-1}^n}; \theta) \, | \, X_{t_{i-1}^n}) = O(\Delta_n^\kappa)$.

\qed

\vspace{5mm}

Theorem \ref{theorem1} follows via asymptotic statistical results for
stochastic processes, see e.g.\ \cite{jjms}. To prove the theorem we
need two technical lemmas. The first is used to establish uniform
convergence in the proofs of Lemma  \ref{lemma2} and Lemma
\ref{lemma3}. The lemma is easier to formulate with the following
definitions.

Let ${\cal C}_0$ denote the subclass of $C_{p,1,2,1}(\R_+,
(\ell,r)^2,\Theta)$ of functions $f(\Delta,y,x;\theta)$ satisfying
that  $f(0,x,x;\theta) =0$ for all $x \in (\ell,r)$ and $\theta \in
\Theta$, and define the operators  
\bean
{\cal L}_1 f(s,y,x;\theta) &=& \partial_s f(s,y,x;\theta) + \partial_y
f(s,y,x;\theta)b(y;\alpha_0) + \mbox{\small $\frac12$} \partial^2_y
f(s,y,x;\theta)v(y,\beta_0) \\
{\cal L}_2 f(s,y,x;\theta) &=& \partial_y f(s,y,x;\theta)\sigma(y;\beta_0).
\eean

\begin{lem}
\label{lemma6}
Assume Condition \ref{con1}, and consider
\be
\label{tightsum}
\zeta^{(j)}_n(\theta) = \frac{1}{n \Delta^{j/2}_n} \sum_{i=1}^n 
f(\Delta_n, X_{t_i^n}, X_{t_{i-1}^n}; \theta), \ \ j= 2,3,4,
\ee
for $f \in {\cal C}_0$. Then the following holds for $j=2$. For every
$m \in \N$ and for every compact $K \subseteq \Theta$, a constant
$C_{m,K} > 0$ exists such that 
\be
\label{momenttightnesscond}
E_{\theta_0} \left( |\zeta_n^{(j)}(\theta_2) - \zeta^{(j)}_n(\theta_1)|^{2m}  
\right) \leq C_{m,K} |\theta_2 - \theta_1|^{2m}
\ee
for all $\theta_1$ and $\theta_2$ in $K$ and for all $n$.

Moreover, if $h_i = {\cal L}_i f \in {\cal C}_0$ for $i=1,2$, then
(\ref{momenttightnesscond}) holds for $j=3$, and if ${\cal
  L}_2 h_i \in {\cal C}_0$ for $i=1,2$, then
(\ref{momenttightnesscond}) holds for $j=4$. 
\end{lem}

\noindent
{\it Proof.} By Ito's formula
\bea
\label{efterIto}
\lefteqn{ f(\Delta_n, X_{t_i^n}, X_{t_{i-1}^n}; \theta) =} \\
&& \int_{t_{i-1}^n}^{t_i^n} h_1(s-t_{i-1}^n, X_{s}, X_{t_{i-1}^n}; \theta)ds +
\int_{t_{i-1}^n}^{t_i^n} h_2(s-t_{i-1}^n, X_{s}, X_{t_{i-1}^n};
\theta)dW_s. \nonumber 
\eea
We treat the two terms on the right hand side of (\ref{efterIto}) separately. 
For a function $k(s,y,x;\theta)$, define $Dk(\cdot ;\theta_2,\theta_1) =
k(\cdot ;\theta_2) - k(\cdot;\theta_1)$. Because $f \in {\cal C}_0$, the partial
derivatives  $\partial_\theta h_i$, $i=1,2$, are of polynomial growth
in $y$ and $x$ uniformly for $\theta$ in a compact set. Therefore
\bean
\lefteqn{\frac{1}{\Delta^{2m}_n} E_{\theta_0} \left( \left| \frac1n 
\sum_{i=1}^n \int_{t_{i-1}^n}^{t_i^n} Dh_1(s-t_{i-1}^n, X_{s}, X_{t_{i-1}^n}; \theta_2,
\theta_1 )ds \right|^{2m} \right)} \\
& \leq & \frac{1}{n \Delta^{2m}_n} \sum_{i=1}^n  E_{\theta_0} 
\left( \left| \int_{t_{i-1}^n}^{t_i^n} Dh_1(s-t_{i-1}^n, X_{s}, X_{t_{i-1}^n}; \theta_2,
\theta_1 )ds \right|^{2m} \right) \\
& \leq & \frac{1}{n \Delta_n} \sum_{i=1}^n  \int_{t_{i-1}^n}^{t_i^n} 
E_{\theta_0}  \left( | Dh_1(s-t_{i-1}^n, X_{s}, X_{t_{i-1}^n}; 
\theta_2, \theta_1 )|^{2m} \right) ds \\
& \leq_C & \frac{1}{n \Delta_n} \sum_{i=1}^n \int_{t_{i-1}^n}^{t_i^n} 
E_{\theta_0} \left( \left | \int_0^1 \partial_\theta h_1(s-t_{i-1}^n, X_{s}, 
X_{t_{i-1}^n}; \theta_1 + u (\theta_2 - \theta_1) )du \right|^{2m}
\right)  ds |\theta_2 - \theta_1|^{2m} \\
& \leq_C & |\theta_2 - \theta_1|^{2m}, 
\eean
where we have used Condition \ref{con1} and Jensen's inequality (twice). 
Using the Burkholder-Davis-Gundy inequality and Jensen's inequality we obtain
\bean
\lefteqn{\frac{1}{\Delta^{2m}_n} E_{\theta_0} \left( \left| \frac1n 
\sum_{i=1}^n \int_{t_{i-1}^n}^{t_i^n} Dh_2(s-t_{i-1}^n, X_{s}, X_{t_{i-1}^n}; \theta_2,
\theta_1 )dW_s \right|^{2m} \right)} \\
& \leq_C &  \frac{1}{\Delta^{2m}_n} E_{\theta_0} \left( \left| \frac{1}{n^2} 
\sum_{i=1}^n \int_{t_{i-1}^n}^{t_i^n} Dh_2(s-t_{i-1}^n, X_{s}, X_{t_{i-1}^n}; \theta_2,
\theta_1)^2 ds \right|^{m} \right) \\
& \leq &  \frac{1}{n^{m+1}\Delta^{2m}_n} \sum_{i=1}^n E_{\theta_0} 
\left( \left| \int_{t_{i-1}^n}^{t_i^n} Dh_2(s-t_{i-1}^n, X_{s}, X_{t_{i-1}^n}; 
\theta_2, \theta_1)^2 ds \right|^{m} \right) \\
& \leq &  \frac{1}{(n \Delta_n)^{m+1}} \sum_{i=1}^n \int_{t_{i-1}^n}^{t_i^n} 
E_{\theta_0} \left( | Dh_2(s-t_{i-1}^n, X_{s}, X_{t_{i-1}^n}; 
\theta_2, \theta_1)|^{2m} \right) ds \\
& \leq_C & \frac{1}{(n \Delta_n)^{m}} |\theta_2 - \theta_1|^{2m},
\eean
which implies (\ref{momenttightnesscond}) for $j=2$.

The result for $j=3,4$ can be proved in a similar way. When $h_i \in
{\cal C}_0$ for $i=1,2$,
\bea
\label{f2}
\lefteqn{f(\Delta_n, X_{t_i^n}, X_{t_{i-1}^n}; \theta) =} \\
&& \int_{t_{i-1}^n}^{t_i^n} \int_{t_{i-1}^n}^{s} 
h_{11}(u-t_{i-1}^n, X_{u}, X_{t_{i-1}^n}; \theta)du ds +
\int_{t_{i-1}^n}^{t_i^n} \int_{t_{i-1}^n}^{s} 
h_{21}(u-t_{i-1}^n, X_{u}, X_{t_{i-1}^n}; \theta)dW_u ds \nonumber \\ 
&& \mbox{} + \int_{t_{i-1}^n}^{t_i^n} \int_{t_{i-1}^n}^{s} 
h_{12}(u-t_{i-1}^n, X_{u}, X_{t_{i-1}^n}; \theta)du dW_s +
\int_{t_{i-1}^n}^{t_i^n} \int_{t_{i-1}^n}^{s} 
h_{22}(u-t_{i-1}^n, X_{u}, X_{t_{i-1}^n}; \theta)dW_u dW_s  \nonumber
\eea
where $h_{ij} = {\cal L}_i {\cal L}_j f$, $i,j = 1,2$. In the two
cases $h_{11}$ and $h_{12}$, we can prove the result for $j=4$, which
implies the result for $j=3$.
\bean
\lefteqn{\frac{1}{\Delta^{4m}_n} E_{\theta_0} \left( \left| \frac1n 
\sum_{i=1}^n \int_{t_{i-1}^n}^{t_i^n}  \int_{t_{i-1}^n}^{s}
Dh_{11}(u-t_{i-1}^n, X_{u}, X_{t_{i-1}^n}; \theta_2, \theta_1 )du ds 
\right|^{2m} \right)} \\
& \leq & \frac{1}{n \Delta^{4m}_n} \sum_{i=1}^n  E_{\theta_0} 
\left( \left| \int_{t_{i-1}^n}^{t_i^n} \int_{t_{i-1}^n}^{s} Dh_{11}(u-t_{i-1}^n, X_{u}, 
X_{t_{i-1}^n}; \theta_2, \theta_1 )du ds \right|^{2m} \right) \\
& \leq & \frac{1}{n \Delta_n^2} \sum_{i=1}^n  \int_{t_{i-1}^n}^{t_i^n} 
\int_{t_{i-1}^n}^{s} E_{\theta_0} \left( | Dh_{11}(u-t_{i-1}^n, X_{u}, 
X_{t_{i-1}^n}; \theta_2, \theta_1 )|^{2m} \right) du ds \\
& \leq_C & |\theta_2 - \theta_1|^{2m}, 
\eean

\bean
\lefteqn{\frac{1}{\Delta^{3m}_n} E_{\theta_0} \left( \left| \frac1n 
\sum_{i=1}^n \int_{t_{i-1}^n}^{t_i^n} \int_{t_{i-1}^n}^{s}
Dh_{21}(u-t_{i-1}^n, X_{u}, X_{t_{i-1}^n}; \theta_2, \theta_1 )dW_u ds 
\right|^{2m} \right)} \\
& \leq &  \frac{1}{\Delta^{3m}_n n}  \sum_{i=1}^n
E_{\theta_0} \left( \left| \int_{t_{i-1}^n}^{t_i^n} \int_{t_{i-1}^n}^{s}
Dh_{21}(u-t_{i-1}^n, X_{u}, X_{t_{i-1}^n}; \theta_2, \theta_1 )dW_u ds 
\right|^{2m} \right) \\
& \leq &  \frac{1}{\Delta^{m+1}_n n}  \sum_{i=1}^n
\int_{t_{i-1}^n}^{t_i^n} E_{\theta_0} \left( \left| \int_{t_{i-1}^n}^{s}
Dh_{21}(u-t_{i-1}^n, X_{u}, X_{t_{i-1}^n}; \theta_2, \theta_1 )dW_u 
\right|^{2m} \right) ds \\
& \leq_C &  \frac{1}{\Delta^{m+1}_n n} \sum_{i=1}^n
\int_{t_{i-1}^n}^{t_i^n} E_{\theta_0} \left( \left| 
\int_{t_{i-1}^n}^{s} Dh_{21}(u-t_{i-1}^n, X_{u}, X_{t_{i-1}^n}; \theta_2,
\theta_1)^2 ds \right|^{m} \right) ds \\
& \leq &  \frac{1}{\Delta_n^{2} n} \sum_{i=1}^n 
\int_{t_{i-1}^n}^{t_i^n} \int_{t_{i-1}^n}^{s} E_{\theta_0} 
\left( | Dh_{21}(u-t_{i-1}^n, X_{u}, X_{t_{i-1}^n}; \theta_2, \theta_1)|^{2m} du 
ds \right) \\
& \leq_C & |\theta_2 - \theta_1|^{2m},
\eean

\bean
\lefteqn{\frac{1}{\Delta^{4m}_n} E_{\theta_0} \left( \left| \frac1n 
\sum_{i=1}^n \int_{t_{i-1}^n}^{t_i^n} \int_{t_{i-1}^n}^{s}
Dh_{12}(u-t_{i-1}^n, X_{u}, X_{t_{i-1}^n}; \theta_2, \theta_1 ) du dW_s 
\right|^{2m} \right)} \\
& \leq_C &  \frac{1}{\Delta^{4m}_n} E_{\theta_0} \left( \left| 
\frac{1}{n^2} \sum_{i=1}^n \int_{t_{i-1}^n}^{t_i^n} \left(
\int_{t_{i-1}^n}^{s} Dh_{12}(u-t_{i-1}^n, X_{u}, X_{t_{i-1}^n}; \theta_2,
\theta_1) du \right)^2 ds \right|^{m} \right) \\
& \leq &  \frac{1}{n^{m+1}\Delta^{4m}_n} \sum_{i=1}^n E_{\theta_0} 
\left( \left| \int_{t_{i-1}^n}^{t_i^n} \left(
\int_{t_{i-1}^n}^{s} Dh_{12}(u-t_{i-1}^n, X_{u}, X_{t_{i-1}^n}; 
\theta_2, \theta_1) du \right)^2 ds \right|^{m} \right) %\\
\eean
\bean
& \leq &  \frac{1}{n^{m+1} \Delta_n^{m+2}} \sum_{i=1}^n 
\int_{t_{i-1}^n}^{t_i^n} \int_{t_{i-1}^n}^{s} E_{\theta_0} \left( 
| Dh_{12}(u-t_{i-1}^n, X_{u}, X_{t_{i-1}^n}; \theta_2, \theta_1)|^{2m} du ds \right) \\
& \leq_C & \frac{1}{(n \Delta_n)^{m}} |\theta_2 - \theta_1|^{2m},
\eean
and

\bean
\lefteqn{\frac{1}{\Delta^{3m}_n} E_{\theta_0} \left( \left| \frac1n 
\sum_{i=1}^n \int_{t_{i-1}^n}^{t_i^n} \int_{t_{i-1}^n}^{s}
Dh_{22}(u-t_{i-1}^n, X_{u}, X_{t_{i-1}^n}; \theta_2, \theta_1 ) dW_u dW_s 
\right|^{2m} \right)} \\
& \leq_C &  \frac{1}{\Delta^{3m}_n} E_{\theta_0} \left( \left| 
\frac{1}{n^2} \sum_{i=1}^n \int_{t_{i-1}^n}^{t_i^n} \left(
\int_{t_{i-1}^n}^{s} Dh_{22}(u-t_{i-1}^n, X_{u}, X_{t_{i-1}^n}; \theta_2,
\theta_1) dW_u \right)^2 ds \right|^{m} \right) \\
& \leq &  \frac{1}{n^{m+1}\Delta^{3m}_n} \sum_{i=1}^n E_{\theta_0} 
\left( \left| \int_{t_{i-1}^n}^{t_i^n} \left(
\int_{t_{i-1}^n}^{s} Dh_{22}(u-t_{i-1}^n, X_{u}, X_{t_{i-1}^n}; 
\theta_2, \theta_1) dW_u \right)^2 ds \right|^{m} \right) \\
& \leq &  \frac{1}{n^{m+1}\Delta^{2m+1}_n} \sum_{i=1}^n 
\int_{t_{i-1}^n}^{t_i^n} E_{\theta_0}  \left( \left|
\int_{t_{i-1}^n}^{s} Dh_{22}(u-t_{i-1}^n, X_{u}, X_{t_{i-1}^n}; 
\theta_2, \theta_1) dW_u \right| ^{2m} \right) ds \\
& \leq_C &  \frac{1}{n^{m+1}\Delta^{2m+1}_n} \sum_{i=1}^n 
\int_{t_{i-1}^n}^{t_i^n} E_{\theta_0}  \left( \left|
\int_{t_{i-1}^n}^{s} Dh_{22}(u-t_{i-1}^n, X_{u}, X_{t_{i-1}^n}; 
\theta_2, \theta_1)^2 du \right| ^{m} \right) ds \\
& \leq &  \frac{1}{n^{m+1} \Delta_n^{m+2}} \sum_{i=1}^n 
\int_{t_{i-1}^n}^{t_i^n} \int_{t_{i-1}^n}^{s} E_{\theta_0} \left( 
| Dh_{22}(u-t_{i-1}^n, X_{u}, X_{t_{i-1}^n}; \theta_2, \theta_1)|^{2m} du ds \right) \\
& \leq_C & \frac{1}{(n \Delta_n)^{m}} |\theta_2 - \theta_1|^{2m}.
\eean

Finally, we prove (\ref{momenttightnesscond}) for $j=4$. We have
already taken care of two of the terms in (\ref{f2}), but the terms
involving $h_{21}$ and $h_{22}$ require more work. Since $h_{2i} \in
{\cal C}_0$, $i=1,2$, we find that
\bean
\int_{t_{i-1}^n}^{t_i^n} \int_{t_{i-1}^n}^{s} 
h_{21}(u-t_{i-1}^n, X_{u}, X_{t_{i-1}^n}; \theta)dW_u ds &=&
\int_{t_{i-1}^n}^{t_i^n} \int_{t_{i-1}^n}^{s} \int_{t_{i-1}^n}^{u}
\hspace{-4mm}
{\cal L}_1h_{21}(v-t_{i-1}^n, X_{v}, X_{t_{i-1}^n}; \theta)dv dW_u ds
\\ &+& \int_{t_{i-1}^n}^{t_i^n} \int_{t_{i-1}^n}^{s} \int_{t_{i-1}^n}^{u}
\hspace{-4mm}
{\cal L}_2 h_{21}(v-t_{i-1}^n, X_{v}, X_{t_{i-1}^n}; \theta)dW_v dW_u ds \\ 
\eean
and
\bean
\int_{t_{i-1}^n}^{t_i^n} \int_{t_{i-1}^n}^{s} 
h_{22}(u-t_{i-1}^n, X_{u}, X_{t_{i-1}^n}; \theta)dW_u dW_s &=&
\int_{t_{i-1}^n}^{t_i^n} \int_{t_{i-1}^n}^{s} \int_{t_{i-1}^n}^{u}
\hspace{-4mm}
{\cal L}_1 h_{22}(v-t_{i-1}^n, X_{v}, X_{t_{i-1}^n}; \theta)dv dW_u
dW_s \\ &+&
\int_{t_{i-1}^n}^{t_i^n} \int_{t_{i-1}^n}^{s} \int_{t_{i-1}^n}^{u}
\hspace{-4mm}
{\cal L}_2 h_{22}(v-t_{i-1}^n, X_{v}, X_{t_{i-1}^n}; \theta)d W_v dW_u dW_s.
\eean
The result is now obtained by evaluating the triple integrals using 
the Burkholder-Davis-Gundy inequality and Jensen's inequality 
exactly as above. 
\qed

\vspace{3mm}

\begin{lem}
\label{lemma2}
Under the Conditions \ref{con1} and \ref{con2}
\be
\label{1lem2}
\frac{1}{n\Delta_n} \sum_{i=1}^n  
g(\Delta_n, X_{t_i^n}, X_{t_{i-1}^n}; \theta) 
\stackrel{P_{\theta_0}}{\longrightarrow} \gamma(\theta,\theta_0),
\ee
\bea
\label{2lem2}
\frac{1}{n\Delta_n} \sum_{i=1}^n \partial_{\theta^T} 
g(\Delta_n, X_{t_i^n}, X_{t_{i-1}^n}; \theta) 
&\stackrel{P_{\theta_0}}{\longrightarrow}& \\
&& \hspace{-6cm} \int_\ell^r 
[L_{\theta_0}(\partial_{\theta^T} g(0;\theta))(x,x) - 
L_{\theta}(\partial_{\theta^T} g(0;\theta))(x,x) - 
J_\theta(x)]\mu_{\theta_0} (x) dx, \nonumber
\eea
and
\be
\label{3lem2}
\frac{1}{n\Delta_n} \sum_{i=1}^n g(\Delta_n, X_{t_i^n}, X_{t_{i-1}^n};
\theta) g(\Delta_n, X_{t_i^n}, X_{t_{i-1}^n}; \theta)^T 
\stackrel{P_{\theta_0}}{\longrightarrow} V(\theta), 
\ee
uniformly for $\theta$ in a compact set. The function $\gamma$ given
by (\ref{gamma}) is a continuous function of $\theta$. For a martingale estimating  
function or more generally if $n \Delta_n^{2\kappa-1} \rightarrow 0$,
\be
\label{4lem2}
\frac{1}{\sqrt{n\Delta_n}} \sum_{i=1}^n 
g(\Delta_n, X_{t_i^n}, X_{t_{i-1}^n}; \theta_0) 
\stackrel{\cal D}{\longrightarrow} 
N_2 \left( 0, V_0 \right).
\ee
\end{lem}

\vspace{4mm}

\noindent
{\it Proof.} By (\ref{gdeltaexp}),
(\ref{expansion}), (\ref{lem1-1}) and Lemma \ref{lemma5},
\bean
\lefteqn{E_{\theta_0} \left( g(\Delta_n, X_{t_i^n}, 
X_{t_{i-1}^n}; \theta) \, | \, X_{t_{i-1}^n} \right)} \\ 
&=& \Delta_n \left[ g^{(1)}(X_{t_{i-1}^n},X_{t_{i-1}^n};
\theta) + L_{\theta_0}(g(0;\theta))(X_{t_{i-1}^n},
X_{t_{i-1}^n}) \right] + \Delta_n^2 R(\Delta_n, X_{t_{i-1}^n},\theta)\\ 
&=& \Delta_n \left[L_{\theta_0}(g(0;\theta))(X_{t_{i-1}^n},
X_{t_{i-1}^n}) - L_{\theta}(g(0;\theta))(X_{t_{i-1}^n},
X_{t_{i-1}^n}) \right] + \Delta_n^2 R(\Delta_n, X_{t_{i-1}^n},\theta).
\eean
The last equality follows from (\ref{lem1-2}). Thus
\bean
\lefteqn{\frac{1}{n\Delta_n} \sum_{i=1}^n  E_{\theta_0} 
\left( g(\Delta_n, X_{t_i^n}, X_{t_{i-1}^n}; \theta) \, | 
\, X_{t_{i-1}^n} \right)} \\
&=& \frac1n \sum_{i=1}^n \left[L_{\theta_0}(g(0;\theta))(X_{t_{i-1}^n},
X_{t_{i-1}^n}) - L_{\theta}(g(0;\theta))(X_{t_{i-1}^n},
X_{t_{i-1}^n}) \right] + \Delta_n \frac1n \sum_{i=1}^n 
R(\Delta_n, X_{t_{i-1}^n},\theta) \\
&& \stackrel{P_{\theta_0}}{\longrightarrow} \gamma(\theta,\theta_0) 
\eean
by Lemma \ref{lemma4}. Moreover, $E_{\theta_0} 
\left( g_j(\Delta_n, X_{t_i^n}, X_{t_{i-1}^n}; \theta)^2 \, | 
\, X_{t_{i-1}^n} \right) = \Delta_n R(\Delta_n, X_{t_{i-1}^n},\theta)$, so
\[
\frac{1}{(n\Delta_n)^2} \sum_{i=1}^n  E_{\theta_0} 
\left( g_j(\Delta_n, X_{t_i^n}, X_{t_{i-1}^n}; \theta)^2 \, | 
\, X_{t_{i-1}^n} \right) = \frac{1}{n\Delta_n} \frac1n
\sum_{i=1}^n R(\Delta_n, X_{t_{i-1}^n},\theta) 
\stackrel{P_{\theta_0}}{\longrightarrow} 0.
\]
Therefore pointwise convergence in (\ref{1lem2}) follows from 
Lemma \ref{lemmagcjj}. In order to prove that the convergence is
uniform for $\theta$ in a compact set $K$, we show that the sequence
$\zeta_n (\cdot) = \frac{1}{n\Delta_n} \sum_{i=1}^n g(\Delta_n,
X_{t_{i}^n},X_{t_{i-1}^n},\cdot)$ converges weakly to the limit 
$\gamma(\cdot,\theta_0)$ in the space, $C(K)$, of continuous functions
on $K$ with the supremum norm. Since the limit is non-random, this implies 
uniform convergence in probability for $\theta \in K$. That
$\gamma(\cdot,\theta_0)$ is continuous follows from the dominated
convergence theorem because of the imposed uniform polynomial growth
assumptions.  Since pointwise convergence has been established, weak
convergence follows because   
the family of distributions of $\zeta_n (\cdot)$ is tight. The tightness
follows from Lemma \ref{lemma6} with $f = g_i$, $j=2$ and $m=2$. 
That (\ref{momenttightnesscond}) and
pointwise convergence implies tightness follows from Corollary 14.9 in
\cite{kallenberg}, which is a generalization of Theorem 12.3 in 
\cite{billingsley} (see also Lemma 3.1 in \cite{nakahiro1} and
Theorem 20 in Appendix I of \cite{ibrahasmin}).

In a similar way it follows from (\ref{gdeltaexp}),
(\ref{expansion}), (\ref{lem1-1}), (\ref{lem1-2}) and 
Lemma \ref{lemma5} that
\bea
\label{sensitivity}
\lefteqn{E_{\theta_0} \left( \partial_{\theta^T} g(\Delta_n, X_{t_i^n}, 
X_{t_{i-1}^n}; \theta) \, | \, X_{t_{i-1}^n} \right)} \\ 
&=& \Delta_n \left[ \partial_{\theta^T} g^{(1)}(X_{t_{i-1}^n},X_{t_{i-1}^n};
\theta) + L_{\theta_0}(\partial_{\theta^T} g(0;\theta))(X_{t_{i-1}^n}, \nonumber
X_{t_{i-1}^n}) \right] + \Delta_n^2 R(\Delta_n, X_{t_{i-1}^n},\theta) \\
&=& \Delta_n \left[L_{\theta_0}(\partial_{\theta^T} g(0;\theta))(X_{t_{i-1}^n},
X_{t_{i-1}^n}) - L_{\theta}(\partial_{\theta^T} g(0;\theta))(X_{t_{i-1}^n},
X_{t_{i-1}^n})  - J_\theta(X_{t_{i-1}^n}) \right] \nonumber \\
&& \hspace{8cm} \mbox{} + \nonumber
\Delta_n^2 R(\Delta_n, X_{t_{i-1}^n},\theta),
\eea 
and from  (\ref{gdeltaexp}),(\ref{expansion}), (\ref{lem1-1}), 
and Lemma \ref{lemma5} that
\bean
\lefteqn{E_{\theta_0} \left( g(\Delta_n, X_{t_i^n}, X_{t_{i-1}^n}; \theta) 
g(\Delta_n, X_{t_i^n}, X_{t_{i-1}^n}; \theta)^T \, 
| \, X_{t_{i-1}^n} \right)} \\ 
&=& \Delta_n v(X_{t_{i-1}^n},\beta_0) 
\partial_y g(0,X_{t_{i-1}^n},X_{t_{i-1}^n};\theta)
\partial_y g(0,X_{t_{i-1}^n},X_{t_{i-1}^n};\theta)^T
+ \Delta_n^2 R(\Delta_n, X_{t_{i-1}^n},\theta).
\eean
Since by (\ref{gdeltaexp}),(\ref{expansion}), (\ref{lem1-1}), 
and Lemma \ref{lemma5}
\[
E_{\theta_0} \left( [\partial_\theta g(\Delta_n, X_{t_i^n}, 
X_{t_{i-1}^n}; \theta)]^2 \, | \, X_{t_{i-1}^n} \right) =
\Delta_n R(\Delta_n, X_{t_{i-1}^n},\theta) 
\]
and
\be
\label{Liapounov1}
E_{\theta_0} \left( [g_j(\Delta_n, X_{t_i^n}, X_{t_{i-1}^n}; \theta) 
g_k(\Delta, X_{t_i^n}, X_{t_{i-1}^n}; \theta)]^2 \, 
| \, X_{t_{i-1}^n} \right)  =
\Delta_n R(\Delta_n, X_{t_{i-1}^n},\theta), 
\ee
we can, as above, use Lemma \ref{lemma4} and  Lemma \ref{lemmagcjj}
to prove (\ref{2lem2}) and (\ref{3lem2}). As above, uniform convergence
for $\theta$ in a compact set $K$ follows by using Lemma \ref{lemma6}
with $f = \partial_{\theta_j} g_k$ and $f = g_j g_k$ to prove the tightness 
of (\ref{tightsum}) (with $j=2$) in $C(K)$.

Finally, (\ref{4lem2}) follows from the central limit theorem for
square integrable martingale arrays under conditions which, in the
martingale case, we have 
already verified in the proof of (\ref{3lem2}), see e.g.\ Corollary 3.1 
in \cite{hallheyde} with the conditional Lindeberg condition replaced
by the stronger conditional Liapounov condition that follows from
(\ref{Liapounov1}) and Lemma \ref{lemma4}, e.g. 
\[
\frac{1}{(n\Delta_n)^2} \sum_{i=1}^n  E_{\theta_0} 
\left( g_j(\Delta_n, X_{t_i^n}, X_{t_{i-1}^n}; \theta_0)^4 \, | 
\, X_{t_{i-1}^n} \right) = \frac{1}{n\Delta_n} \frac1n
\sum_{i=1}^n R(\Delta_n, X_{t_{i-1}^n},\theta_0) 
\stackrel{P_{\theta_0}}{\longrightarrow} 0.
\]
The nestedness condition in Hall and Heyde's Corollary 3.1 is not 
needed here because the limit of the quadratic variation
is non-random.

In the case of non-martingale estimating functions, we
consider the martingale $\sum_{i=1}^n \tilde g(\Delta_n,$ $X_{t_i^n}, 
 X_{t_{i-1}^n}; \theta_0)$, where $\tilde g = g - E_{\theta_0} 
\left( g | \, X_{t_{i-1}^n} \right)$. This martingale satisfies the conditions of 
the central limit theorem, which follows from the expansions of
conditional expectations given above and  $E_{\theta_0} 
\left( g_j(\Delta_n, X_{t_i^n}, X_{t_{i-1}^n}; \theta_0)^3 \, | 
\, X_{t_{i-1}^n} \right) = \Delta_n R(\Delta_n, X_{t_{i-1}^n},\theta_0)$.
Now,  (\ref{4lem2}) follows because by (\ref{mart})
\pagebreak
\be
\label{non-mart-cond}
\frac{1}{\sqrt{n \Delta_n}} \sum_{i=1}^n E_{\theta_0} \left(
g(\Delta_n, X_{t_i^n}, X_{t_{i-1}^n}; \theta_0)  \, 
| \, X_{t_{i-1}^n} \right)  =  \sqrt{n} \Delta_n^{\kappa-1/2}
\frac1n \sum_{i=1}^n R(\Delta_n, X_{t_{i-1}^n},\theta_0) 
\stackrel{P_{\theta_0}}{\longrightarrow} 0.
\ee
\qed

\vspace{3mm}

\noindent
{\it Proof of Theorem \ref{theorem1}.} By Lemma \ref{lemma2}, the 
estimating function
\be
\label{G}
G_n(\theta) = \frac{1}{n\Delta_n} \sum_{i=1}^n  
g(\Delta_n, X_{t_i^n}, X_{t_{i-1}^n}; \theta) 
\ee
satisfies the conditions that $G_n(\theta_0) 
\stackrel{P_{\theta_0}}{\rightarrow} 0$, $\partial_\theta G_n(\theta) 
\stackrel{P_{\theta_0}}{\rightarrow} U(\theta)$ uniformly for $\theta$
in a compact set, and that $U(\theta_0) = -S$ is invertible. Here 
$U(\theta)$ denotes the right hand side of (\ref{2lem2}). This implies
the eventual existence and the consistency of $\hat{\theta}_n$ as well
as the eventual uniqueness of consistent $G_n$-estimators; see
Theorems 2.5 and 2.6 in \cite{jjms}. Now consider any $G_n$-estimator
$\tilde \theta_n$ for which $P_{\theta_0}(\tilde \theta_n \in K)
\rightarrow 1$ as $n  \rightarrow \infty$, where $K$ is a compact
subset of $\Theta$ with $\theta_0 \in$ int$\, K$. By Theorem
2.7 in \cite{jjms} the facts that $\gamma(\theta, \theta_0)$ (the
limit of $G_n(\theta)$) satisfies that $\gamma(\theta,\theta_0) \neq 0$
for $\theta \neq \theta_0$ and is continuous in $\theta$ imply that
$P_{\theta_0}(\hat \theta_n \neq \tilde \theta_n) \rightarrow 0$ as $n
\rightarrow \infty$. The asymptotic normality follows by standard
arguments. Finally, (\ref{varest1}) and (\ref{varest2}) follow from
(\ref{2lem2}) and (\ref{3lem2}) because the convergence is uniform for
$\theta$ in compact sets. \qed

\vspace{3mm}

The next lemma is needed in the proof of Theorem \ref{theorem2}

\begin{lem}
\label{lemma3}
Under the Conditions \ref{con3}, \ref{con1}, and \ref{con2}
\be
\label{1lem3}
D_n \sum_{i=1}^n g(\Delta_n, X_{t_i^n}, X_{t_{i-1}^n};
\theta) g(\Delta_n, X_{t_i^n}, X_{t_{i-1}^n}; \theta)^T D_n 
\stackrel{P_{\theta_0}}{\longrightarrow}
\left( \begin{array}{cc} W_1(\theta) & 0 \\ 0 & W_2(\theta) \rule{0mm}{6mm}
\end{array} \right)
\ee
uniformly for $\theta$ in a compact set, where $D_n$ is given by
(\ref{Dn}). \\

For a martingale estimating function or if more generally 
$n \Delta^{2(\kappa-1)} \rightarrow 0$,
\be
\label{2lem3}
\left( \begin{array}{c} \frac{1}{\sqrt{n\Delta_n}} 
\sum_{i=1}^n g_1(\Delta_n, X_{t_i^n}, X_{t_{i-1}^n}; \theta_0) \\
\frac{1}{\Delta_n \sqrt{n}} 
\sum_{i=1}^n g_2(\Delta_n, X_{t_i^n}, X_{t_{i-1}^n}; 
\theta_0) \rule{0mm}{7mm}
\end{array} \right)
\stackrel{\cal D}{\longrightarrow} 
N_2 \left( \left( \begin{array}{c} 0 \\ 0 \end{array} \right), 
\left( \begin{array}{cc} W_1(\theta_0) & 0 \\ 0 & 
W_2(\theta_0) \rule{0mm}{6mm} \end{array} \right) \right).
\ee
If, in addition, condition (\ref{extracond}) holds, then
\be
\label{extraconv}
\frac{1}{n\Delta^{3/2}_n} \sum_{i=1}^n \partial_{\alpha} 
g_2(\Delta_n, X_{t_i^n}, X_{t_{i-1}^n}; \theta) 
\stackrel{P_{\theta_0}}{\longrightarrow} 0
\ee
uniformly for $\theta$ in a compact set.

\end{lem}

\noindent
{\it Proof.} 
By (\ref{gdeltaexp}), (\ref{expansion}), (\ref{lem1-1}), 
(\ref{jacobsen}) and Lemma \ref{lemma5},
\bean
\lefteqn{\frac{1}{n\Delta_n^{3/2}} \sum_{i=1}^n  E_{\theta_0} 
\left( g_1(\Delta_n, X_{t_i^n}, X_{t_{i-1}^n}; \theta)
g_2(\Delta_n, X_{t_i^n}, X_{t_{i-1}^n}; \theta) \, | 
\, X_{t_{i-1}^n} \right)} \\
&&  \hspace{6cm}= \Delta_n^{1/2} \frac1n
\sum_{i=1}^n R(\Delta_n, X_{t_{i-1}^n},\theta) 
\stackrel{P_{\theta_0}}{\longrightarrow} 0
\eean
and
\bea
\label{Liapounov3}
\lefteqn{\frac{1}{n^2\Delta_n^3} \sum_{i=1}^n  E_{\theta_0} 
\left( [g_1(\Delta_n, X_{t_i^n}, X_{t_{i-1}^n}; \theta)
g_2(\Delta_n, X_{t_i^n}, X_{t_{i-1}^n}; \theta)]^2 \, | 
\, X_{t_{i-1}^n} \right)} \\
&&  \hspace{6cm}= \frac{1}{n\Delta_n} \frac1n
\sum_{i=1}^n R(\Delta_n, X_{t_{i-1}^n},\theta) 
\stackrel{P_{\theta_0}}{\longrightarrow} 0, \nonumber
\eea
so the pointwise convergence of the two off-diagonal entries in (\ref{1lem3})
follows from Lemma \ref{lemmagcjj}.
Similarly to the proof of Lemma \ref{lemma2}, uniform convergence
for $\theta$ in a compact set $K$ follows by using Lemma \ref{lemma6}
with $f = g_1 g_2$ to prove the tightness of (\ref{tightsum})
(with $j=3$) in $C(K)$.

The convergence of 
$(n \Delta_n)^{-1} \sum_{i=1}^n g_1(\Delta_n, X_{t_i^n}, X_{t_{i-1}^n};
\theta)^2 $ was taken care of in Lemma \ref{lemma2}. By (\ref{gdeltaexp}),
(\ref{expansion}), (\ref{lem1-1}), (\ref{lem1-2}), 
(\ref{jacobsen}) and Lemma \ref{lemma5},
we see that
\bean
\lefteqn{E_{\theta_0} \left( g_2(\Delta_n, X_{t_i^n}, 
X_{t_{i-1}^n}; \theta)^2 \, | \, X_{t_{i-1}^n} \right)} \\ 
&=& \Delta_n^2 \left[  \mbox{\small $\frac12$} L^2_{\theta_0}
(g_2(0;\theta)^2)(X_{t_{i-1}^n}, X_{t_{i-1}^n}) + 2 L_{\theta_0}
(g_2(0;\theta)g_2^{(1)}(\theta))(X_{t_{i-1}^n}, X_{t_{i-1}^n})
\right. \\ && \hspace{6cm} \left. \mbox{}
+ g_2^{(1)}(X_{t_{i-1}^n}, X_{t_{i-1}^n};\theta)^2
\right] + \Delta_n^3 R(\Delta_n, X_{t_{i-1}^n},\theta) \\
&=& \mbox{\small $\frac12$} \Delta_n^2 \left[ v(X_{t_{i-1}^n};\beta_0)^2
+ \mbox{\small $\frac12$} (v(X_{t_{i-1}^n};\beta_0) -
v(X_{t_{i-1}^n};\beta))^2 \right] (\partial_y^2 g_2(0, X_{t_{i-1}^n}, 
X_{t_{i-1}^n};\theta))^2 \\
&& \hspace{10cm} \mbox{} + \Delta_n^3 R(\Delta_n, X_{t_{i-1}^n},\theta),
\eean 
Thus
\bean
\lefteqn{\frac{1}{n\Delta_n^2} \sum_{i=1}^n  E_{\theta_0} 
\left( g_2(\Delta_n, X_{t_i^n}, X_{t_{i-1}^n}; \theta)^2 \, | 
\, X_{t_{i-1}^n} \right)} \\
&=& \frac1n \sum_{i=1}^n \mbox{\small $\frac12$} 
\left[ v(X_{t_{i-1}^n};\beta_0)
+ \mbox{\small $\frac12$} (v(X_{t_{i-1}^n};\beta_0) -
v(X_{t_{i-1}^n};\beta))^2 \right] (\partial_y^2 g_2(0, X_{t_{i-1}^n}, 
X_{t_{i-1}^n};\theta))^2 \\ && \hspace{8cm}
\mbox{} + \Delta_n \frac1n \sum_{i=1}^n 
R(\Delta_n, X_{t_{i-1}^n},\theta) \\
&& \stackrel{P_{\theta_0}}{\longrightarrow} W_2(\theta) 
\eean
by Lemma \ref{lemma4}. We conclude that $(n \Delta_n^2)^{-1} 
\sum_{i=1}^n g_2(\Delta_n, X_{t_i^n}, X_{t_{i-1}^n} \theta)^2 $  
converges to $W_2(\theta)$ by Lemma \ref{lemmagcjj} because
\be
\label{Liapounov2}
\frac{1}{n^2\Delta_n^4} \sum_{i=1}^n  E_{\theta_0} 
\left( g_2(\Delta_n, X_{t_i^n}, X_{t_{i-1}^n}; \theta)^4 \, | 
\, X_{t_{i-1}^n} \right) = \frac{1}{n\Delta_n} \frac1n
\sum_{i=1}^n R(\Delta_n, X_{t_{i-1}^n},\theta) 
\stackrel{P_{\theta_0}}{\longrightarrow} 0.
\ee
This follows from (\ref{gdeltaexp}), (\ref{expansion}), (\ref{lem1-1}),  
(\ref{jacobsen}), and Lemmas \ref{lemma5} and \ref{lemma4}.
Uniform convergence for $\theta$ in a compact set $K$ 
follows by using Lemma \ref{lemma6} with $f = g_2^2$ to prove 
the tightness of (\ref{tightsum}) (with $j=4$) in $C(K)$.

As in the proof of Lemma \ref{lemma2}, (\ref{2lem3}) follows from 
the central limit theorem for square integrable 
martingale arrays (Corollary 3.1 in \cite{hallheyde}) under conditions 
which, in the martingale case, we have 
already verified in the proof of (\ref{1lem3}). In particular, the
conditional Liapounov condition follows from (\ref{Liapounov1}),
(\ref{Liapounov2}) and (\ref{Liapounov3}).

In the case of non-martingale estimating functions, consider the martingale
$\sum_{i=1}^n \tilde g(\Delta_n,$ $ X_{t_i^n},  X_{t_{i-1}^n};
\theta_0)$, where $\tilde g = g - E_{\theta_0}  \left( g | \,
  X_{t_{i-1}^n} \right)$, which satisfies the conditions of  the central limit theorem.
This follows from the expansions of conditional expectations given
above and  $E_{\theta_0}  \left( g_2(\Delta_n, X_{t_i^n}, X_{t_{i-1}^n}; \theta_0)^3 \, | 
\, X_{t_{i-1}^n} \right) = \Delta_n^2 R(\Delta_n, X_{t_{i-1}^n},\theta_0)$.
Now, (\ref{2lem3}) follows because $g_1$
satisfies (\ref{non-mart-cond}), and because by (\ref{mart})
\[
\frac{1}{\Delta_n\sqrt{n}} \sum_{i=1}^n E_{\theta_0} \left(
g_2(\Delta_n, X_{t_i^n}, X_{t_{i-1}^n}; \theta_0)  \, 
| \, X_{t_{i-1}^n} \right)  =  \sqrt{n} \Delta_n^{\kappa-1}
\frac1n \sum_{i=1}^n R(\Delta_n, X_{t_{i-1}^n},\theta_0) 
\stackrel{P_{\theta_0}}{\longrightarrow} 0.
\]

Finally, to prove (\ref{extraconv}) note that (\ref{sensitivity}),
(\ref{jacobsen}) and (\ref{extracond}) imply that
\[
E_{\theta_0} \left( \partial_\alpha g_2(\Delta_n, X_{t_i^n}, 
X_{t_{i-1}^n}; \theta) \, | \, X_{t_{i-1}^n} \right)
= \Delta^2 R(\Delta_n, X_{t_{i-1}^n},\theta), 
\]
and that it follows from (\ref{gdeltaexp}),
(\ref{expansion}), (\ref{lem1-1}), (\ref{lem1-2}), 
(\ref{jacobsen}), (\ref{extracond}) and Lemma \ref{lemma5} that
\[
E_{\theta_0} \left( [\partial_\alpha g_2(\Delta_n, X_{t_i^n}, 
X_{t_{i-1}^n}; \theta)]^2 \, | \, X_{t_{i-1}^n} \right)
= \Delta_n^3 R(\Delta_n, X_{t_{i-1}^n},\theta). 
\]
Therefore by Lemma \ref{lemma4}
\[
\frac{1}{n \Delta_n^{3/2}} \sum_{i=1}^n  E_{\theta_0} 
\left( \partial_\alpha g_2(\Delta_n, X_{t_i^n}, X_{t_{i-1}^n}; 
\theta) \, | \, X_{t_{i-1}^n} \right) =
\sqrt{\Delta_n}\frac1n \sum_{i=1}^n R(\Delta_n,X_{t_{i-1}^n};\theta)
\stackrel{P_{\theta_0}}{\longrightarrow} 0.
\]
and
\[
\frac{1}{n^2 \Delta_n^3} \sum_{i=1}^n  E_{\theta_0} 
\left( [\partial_\alpha g_2(\Delta_n, X_{t_i^n}, X_{t_{i-1}^n}; 
\theta)]^2 \, | \, X_{t_{i-1}^n} \right) =
\frac{1}{n} \frac1n \sum_{i=1}^n 
R(\Delta_n, X_{t_{i-1}^n},\theta) 
\stackrel{P_{\theta_0}}{\longrightarrow} 0,
\]
so that (\ref{extraconv}) follows from Lemma \ref{lemmagcjj}. 
Uniform convergence for $\theta$ in a compact set $K$ follows by
using Lemma \ref{lemma6} with $f = \partial_\alpha g_2$ to conclude 
tightness of (\ref{tightsum}) (with $j=3$) in $C(K)$. To see that 
$\partial_\alpha g_2$ satisfies the conditions of the lemma, we use 
(\ref{lem1-2}) and (\ref{extracond}) to conclude that $\partial_\Delta \partial_\alpha 
g_2(0,x,x;\theta) = \partial_\alpha g_2^{(1)}(x,x;\theta) =
- \partial_\alpha L_\theta(g_2(0;\theta))(x,x)=0$.

\qed

\vspace{3mm}

\noindent
{\it Proof of Theorem \ref{theorem2}.} The results on eventual existence,
uniqueness and consistence of $\hat{\theta}_n$ follow from Theorem \ref{theorem1}:
Because (\ref{jacobsen}) implies $S_{21} = 0$, the assumptions that 
$S_{11} \neq 0$ and $S_{22} \neq 0$ ensure that $S$ is invertible,
and similarly, under Condition \ref{con3} the identifiability condition imposed
in Theorem \ref{theorem2} ensures that $\gamma(\theta,\theta_0) 
\neq 0$ for $\theta \neq
\theta_0$, where $\gamma$ is the limit of $G_n(\theta)$ given by
(\ref{gamma}).

To prove the asymptotic normality (\ref{asnorm2}) of the estimator 
$\hat{\theta}_n$ we consider
\[
\tilde{G}_n(\theta) = D_n \sum_{i=1}^n  g(\Delta_n, X_{t_i^n}, 
X_{t_{i-1}^n}; \theta), 
\]
where $D_n$ is given by (\ref{Dn}). On the set $\{ 
\tilde{G}_n(\hat{\theta}_n) = 0 \}$ (the probability of which goes to
one)
\[
- \partial_{\theta^T} \tilde{G}_n(\theta^{(1)}_n,\theta^{(2)}_n)
A_n^{-1}A_n (\hat{\theta}_n - \theta_0) = \tilde{G}_n(\theta_0),
\]
where
\[
A_n = \left( \begin{array}{cc} \sqrt{\Delta_n n} & 0 \\ 0 &
\sqrt{n} \end{array} \right),
\]
$\partial_{\theta^T} \tilde{G}_n(\theta^{(1)}_n,
\theta^{(2)}_n)$ is the $2 \times 2$-matrix whose $jk$th entry is
$\partial_{\theta_k} \tilde{G}_n(\theta^{(j)}_n)_j$, and
$\theta^{(j)}_n$ is a random convex combination of $\hat{\theta}_n$
and $\theta_0$. Since by (\ref{2lem2}) and (\ref{extraconv})
\[
- \partial_{\theta^T} \tilde{G}_n(\theta^{(1)}_n,\theta^{(2)}_n)
A_n^{-1} \stackrel{P_{\theta_0}}{\longrightarrow}
\left( \begin{array}{cc} S_{11} & 0 \\ 0 & S_{22} \end{array} \right),
\]
(\ref{asnorm2}) follows from (\ref{2lem3}).

Finally, (\ref{varest3}) follows from (\ref{1lem3}) because the
convergence is uniform for $\theta$ in compact sets.

\qed 

\vspace{3mm}

\noindent
{\it Proof of Theorem \ref{theorem:efficiency}.} That the conclusions
of Theorem \ref{theorem2} holds is obvious because (\ref{diffeff})
implies (\ref{extracond}). The efficiency follows from Theorem 4.1 in
\cite{gobet}, where it is proved that the diffusion model (\ref{sde}) is locally 
asymptotically normal with Fisher information matrix ${\cal
  I}(\theta_0) = \Sigma(\theta_0)^{-1}$, with $\Sigma(\theta_0)$ given
by (\ref{optimalcovariance}). Under Condition \ref{conefficiency},
$S_{11} = W_1(\theta_0) = {\cal I}_{11}(\theta_0)$ and  $S_{22} = W_2(\theta_0)
= {\cal I}_{22}(\theta_0)$, so it follows from (\ref{asnorm2}) that the
asymptotic covariance matrix of $\hat \theta_n$ equals the inverse of the
Fisher information matrix. The estimators of the asymptotic variances follow
from (\ref{varest3}).

\qed

\vspace{3mm}

\noindent
{\it Proof of Theorem \ref{theorem4}.} By (\ref{expansion})
$\pi^{1,\Delta}_\theta f(x;\theta) = f(x;\theta) + \Delta L_\theta 
f(x;\theta) + O(\Delta^2)$,
so after another application of (\ref{expansion}), we see that
$h(\Delta,y,x;\theta) = f(y;\theta) - \pi_\theta^{\kappa,\Delta} 
f(x;\theta)$, $\kappa \in \N$, satisfies
\bean
%\label{Ehh}
\lefteqn{ E_\theta \left( h(\Delta, X_\Delta, x; \theta)
h(\Delta, X_\Delta, x; \theta)^T \, | \, X_0 = x \right)
= \Delta L_\theta ( h(0;\theta)h(0;\theta)^T)(x,x)} \\ 
&& \mbox{} \hspace{10mm} + \Delta^2
\left( \mbox{\small $\frac12$} L^2_\theta ( h(0;\theta)
h(0;\theta)^T)(x,x) - L_\theta f(x;\theta) L_\theta f^T(x;\theta)
\right) + O(\Delta^3) \\
&=& \Delta v(x;\beta) \partial_x f(x;\theta)\partial_x f(x;\theta)^T
+ \Delta^2 K(x) + O(\Delta^3),
\eean
where
\bean
K(x) &=& q_1(x;\theta)\partial_x f(x;\theta)\partial_x f(x;\theta)^T 
+ q_2(x;\theta) \left(\partial^2_x f(x;\theta)\partial_x f(x;\theta)^T 
+ \partial_x f(x;\theta)\partial^2_x f(x;\theta)^T \right) \\ 
&& \mbox{} %\hspace{-7.5mm} 
+ v(x;\beta)^2 
\left( \partial^2_x f(x;\theta)\partial^2_x f(x;\theta)^T + 
\mbox{\small $\frac12$} (\partial^3_x f(x;\theta)\partial_x f(x;\theta)^T 
+ \partial_x f(x;\theta)\partial^3_x f(x;\theta)^T) \right) \nonumber
\eean
with
\vspace{-2mm}
\bean
q_1(x;\theta) &=& \mbox{\small $\frac12$}[b(x;\alpha)(2 + \partial_x
v(x;\beta)) - 2 b (x;\alpha)+  \mbox{\small $\frac12$} v(x;\beta)
(4 \partial_x b(x;\alpha) + \partial_x^2 v(x;\beta))] \\
q_2(x,\theta) &=& \mbox{\small $\frac34$} v(x;\beta)(1 + 
\mbox{\small $\frac13$}b(x;\alpha) + \partial_x v(x;\beta)).
\eean
Since
\vspace{-3mm}
\bean
\partial_\alpha L_\theta f(x;\theta) - L_\theta \partial_\alpha f(x;\theta) 
&=& \partial_\alpha b(x;\alpha) \partial_x f(x;\theta) \\
\partial_\beta L_\theta f(x;\theta) -  L_\theta \partial_\beta f(x;\theta)
&=& \mbox{\small $\frac12$} 
\partial_\beta v(x;\beta) \partial^2_xf(x;\theta),
\eean
it also follows from (\ref{expansion}) that
\[
\partial_{\theta^T} \pi^\Delta_\theta f(x) -
\pi^\Delta_\theta \partial_{\theta^T} f(x) = \Delta F(x)  
\left( \begin{array}{cc} \partial_\alpha b(x;\alpha) & 0 \\
0 & \frac12 \partial_\beta v(x;\beta)\rule{0mm}{6mm} \end{array} \right) + 
O(\Delta^2),
\]
where $F(x)$ denotes the $N \times 2$-matrix $F(x) = ( \partial_xf(x),
\partial^2_xf(x) )$.

If $A^*(x,\Delta;\theta)$ satisfies (\ref{godambeheyde}), then
the $2 \times N$-matrix
\[
B(x,\Delta;\theta) = \left( \begin{array}{cc}
1 & 0 \\ 0 & 2 \Delta \rule{0mm}{5mm} \end{array} \right)  
A^* (x, \Delta; \theta).
\]
satisfies that
\bea
\label{eq2}
\lefteqn{B(x,\Delta;\theta)\left[ v(x;\beta)\partial_x f(x;\theta)
\partial_x f(x;\theta)^T + \Delta K(x;\theta) + O(\Delta^2)\right]} \\ 
&& \nonumber \\
&& \hspace{3cm} = 
\left( \begin{array}{cc} \partial_\alpha b(x;\alpha) & 0 \\
0 & \Delta \partial_\beta v(x;\beta) \rule{0mm}{6mm}
\end{array} \right) F(x)^T 
+ \left( \begin{array}{c} O(\Delta) \\ O(\Delta^2) 
\rule{0mm}{6mm} \end{array} \right). \nonumber
\eea
Let $B(x,\Delta;\theta)_i$ denote the $i$th row of
$B(x,\Delta;\theta)$ ($i=1,2$). Then it follows by letting $\Delta$
tend to zero that
\be
\label{eq1}
v(x;\beta) B(x,0;\theta)_2 \partial_x f(x;\theta)\partial_x f(x;\theta)^T = 0.
\ee
The condition that $M(x)$ is invertible implies that we can find a 
coordinate of $\partial_x f(x;\theta)$ which is not equal to zero, 
so we conclude that
\[
\partial_y g_2^*(0,x,x;\theta) = B(x,0;\theta)_2 \partial_x f(x;\theta) = 0.
\]
Similarly we find that
\[
[v(x;\beta) B(x,0;\theta)_1 \partial_x f(x;\theta) - \partial_\alpha 
b(x;\alpha)] \partial_x f(x;\theta)^T = 0,
\]
which implies
\[
\partial_y g_1^*(0,x,x;\theta) = B(x,0;\theta)_1 \partial_x f(x;\theta) = 
\partial_\alpha b(x;\alpha)/v(x;\beta).
\]
Finally, (\ref{eq2}) and (\ref{eq1}) imply that $B(x,\Delta;\theta)_2
\left(\Delta K(x;\theta) + O(\Delta^2)\right)= \Delta \partial_\beta
v(x;\beta)\partial^2_x f(x;\theta)^T+ O(\Delta^2)$, so that
\[
B(x,0;\theta)_2 K(x;\theta) = \partial_\beta v(x;\beta)\partial^2_x 
f(x;\theta)^T. 
\]
Since we have shown that $B(x,0;\theta)_2 \partial_x f(x;\theta) = 0$,
this expression can be rewritten as
\[
c_1(x;\theta)\partial_x f(x;\theta) = c_2(x;\theta)\partial^2_x 
f(x;\theta)
\]
where

\vspace{-8mm}

\bean
c_1(x;\theta) &=& q_2(x;\theta)B(x,0;\theta)_2 \partial^2_x f(x;\theta) +
\mbox{\small $\frac12$} v(x;\beta)^2 B(x,0;\theta)_2 \partial^3_x
f(x;\theta) \\
c_2(x;\theta) &=& \partial_\beta v(x;\beta) - v(x;\beta)^2 
B(x,0;\theta)_2 \partial^2_x f(x;\theta).
\eean
If $c_2(x;\theta) \neq 0$, then $\partial^2_x f(x) =
c_1(x;\theta)/c_2(x;\theta) \partial_x f(x;\theta)$, which implies 
that det$ (M(x))=0$. This contradicts the assumption that $M(x)$ is
invertible, so we conclude that $\partial_\beta 
v(x;\beta) - v(x;\beta)^2  B(x,0;\theta)_2 \partial^2_x f(x;\theta) 
= 0$ or equivalently
\[
\partial^2_y g_2^*(0,x,x;\theta) = B(x,0;\theta)_2 \partial^2_x f(x;\theta) = 
\partial_\beta v(x;\beta)/v(x;\beta)^2. 
\]

That the results hold for $\tilde B$ is obvious. Note that it also
follows that neither $B(x,0;\theta)_1$ nor $B(x,0;\theta)_2$ is the
zero vector, so there exist $x$ and $y$ such that $g_i^*(0,y,x;\theta)
\neq 0$, $i=1,2$.

\qed

\section*{Acknowledgements}

I am grateful to Arnaud Gloter, Nina Munkholt Jakobsen, Mathieu
Kessler, Masayuki Uchida and Nakahiro Yoshida for helpful discussions. 

\bibliographystyle{natbib}
\bibliography{efficient}

\end{document}